\documentclass[12pt,a4paper]{article}
\usepackage{amsmath,amsfonts,amssymb,amsthm,amscd}
\usepackage{hyperref}
\usepackage{braket}
\usepackage{epstopdf}

\usepackage{braket}
\usepackage{xypic}
\usepackage{amsmath}

\usepackage[american]{babel}
\usepackage{color} 
\usepackage{epsfig}
\usepackage{caption}
\usepackage{rotating}
\usepackage{setspace}
\usepackage{fancyhdr}
\usepackage{booktabs}
\usepackage{mathrsfs}
\usepackage{amsthm}
\usepackage{amsmath,amssymb}
\usepackage{enumerate}
\usepackage{tikz}
\usepackage{youngtab}
\usepackage{hyperref}
\usepackage{tensor}
\usepackage{mathabx}
\usepackage{booktabs}
\usepackage{mathrsfs}
\usepackage{subfigure}
\usepackage{amsmath,amssymb}

\renewcommand{\[}{\begin{equation}}
\renewcommand{\]}{\end{equation}}

\newtheorem{thm}{Theorem}[section]
\newtheorem{cor}[thm]{Corollary}
\newtheorem{lem}[thm]{Lemma}
\newtheorem{prop}[thm]{Proposition}

\theoremstyle{definition}
\newtheorem{defn}{Definition}[section]

\theoremstyle{remark}

\DeclareMathOperator{\PSL}{PSL}

\newcommand{\ra}{\rightarrow}

\newcommand{\zz}{\mathbb{Z}}
\newcommand{\nn}{\mathbb{N}}
\newcommand{\cc}{\mathbb{C}}
\newcommand{\rr}{\mathbb{R}}

\DeclareMathOperator{\ran}{ran}
\DeclareMathOperator{\Tr}{Tr}

\DeclareMathOperator{\SL}{SL}
\DeclareMathOperator{\PGL}{PGL}
\DeclareMathOperator{\SO}{SO}
\DeclareMathOperator{\GL}{GL}

\DeclareMathOperator{\rank}{rank}
\DeclareMathOperator{\gz}{G(\mathbb{Z})}

\DeclareMathOperator{\sign}{sign}
\DeclareMathOperator{\inv}{inv}
\DeclareMathOperator{\Inv}{Inv}
\DeclareMathOperator{\topp}{top}

\DeclareMathOperator{\g}{\gamma^{(0)}}
\DeclareMathOperator{\guno}{\gamma^{(1)}}
\DeclareMathOperator{\gdue}{\gamma^{(2)}}
\DeclareMathOperator{\e}{\eta^{(0)}}
\DeclareMathOperator{\euno}{\eta^{(1)}}
\DeclareMathOperator{\edue}{\eta^{(2)}}

\DeclareMathOperator{\lin}{lin}

\title{A mixing property for the action of $\SL(3,\mathbb{Z}) \times \SL(3,\mathbb{Z})$ on the Stone-{\v C}ech boundary of $\SL(3,\mathbb{Z})$}

\author{J. Bassi and F. R{\u a}dulescu\footnote{Florin R{\u a}dulescu is a member of the Institute of Mathematics of the Romanian Academy} \\ \\Department of Mathematics, University of Tor Vergata, \\Via della Ricerca Scientifica 1, 00133, Roma, Italy\\ \\email: bssjcp01@uniroma2.it; radulesc@mat.uniroma2.it}

\begin{document}

\date{}
\maketitle

\begin{abstract}
\noindent By analogy with the construction of the Furstenberg boundary, the Stone-{\v C}ech boundary of $\SL(3,\mathbb{Z})$ is a fibered space over products of projective matrices. The proximal behaviour on this space is exploited to show that the preimages of certain sequences have accumulation points which belong to specific regions, defined in terms of flags. We show that the $\SL(3,\mathbb{Z})\times \SL(3,\mathbb{Z})$-quasi-invariant Radon measures supported on these regions are tempered. Thus every quasi-invariant Radon boundary measure for $\SL(3,\mathbb{Z})$ is an orthogonal sum of a tempered measure and a measure having matrix coefficients belonging to a certain ideal $c'_0 ((\SL(3,\mathbb{Z}) \times \SL(3,\mathbb{Z}))$, slightly larger than $c_0 ((\SL(3,\mathbb{Z}) \times \SL(3,\mathbb{Z}))$. Hence the left-right representation of $C^*(\SL(3,\mathbb{Z}) \times \SL(3,\mathbb{Z}))$ in the Calkin algebra of $\SL(3,\mathbb{Z})$ factors through $C^*_{c'_0} (\SL(3,\mathbb{Z}) \times \SL(3,\mathbb{Z}))$ and the centralizer of every infinite subgroup of $\SL(3,\mathbb{Z})$ is amenable. 





\end{abstract}
\tableofcontents

\section{Introduction}
 In their seminal paper (\cite{AkOs}) Akemann and Ostrand observed an interesting fact concerning the free group $\mathbb{F}_2$: while the left-right representation of $C^* (\mathbb{F}_2 \times \mathbb{F}_2)$ on $l^2 (\mathbb{F}_2)$ cannot factor through a representation of $C^*_{red} (\mathbb{F}_2 \times \mathbb{F}_2)$ (this can only happen for amenable groups), this factorization appears after passing to the Calkin algebra. Groups for which this phenomenon occurs are said to satisfy the Akemann-Ostrand property and the associated von Neumann algebras are solid (\cite{ozas}). Examples are given by hyperbolic groups and counterexamples by groups containing infinite subgroups with non-amenable centralizers (\cite{Sk}). In particular $\SL(n,\zz)$ satisfies the Akemann-Ostrand property for $n=2$ and does not for $n \geq 4$, while the case $n=3$ is still an open problem to the authors' knowledge. This work is aimed to a better understanding of the topology of the image of $C^*(\SL(3,\zz) \times \SL(3,\zz))$ in the Calkin algebra $\mathbb{B}(l^2 (\SL(3,\zz)))/\mathbb{K}(l^2 (\SL(3,\zz)))$. In contrast with $\SL(2,\mathbb{Z})$, the groups $\SL(n,\mathbb{Z})$, with $n \geq 3$, all have property $T$; moreover they also share some important rigidity properties (see for instance \cite{Zi} Chapter 5, \cite{F} and \cite{F1}). The results contained in this work show that the case $n=3$ also has some common feature with the case $n=2$, which are not shared in higher rank.
  
 The description of the Stone-{\v C}ech boundary as a fibered space developed here is inspired by some considerations of Furstenberg appearing in \cite{Fu}. In order to motivate the discussion, we briefly review the case of $\SL(2,\zz)$ within Furstenberg's approach (c.f. \cite{Fu} pp. 34--35). Let $\gamma_n \rightarrow \infty$ in $\SL(2,\zz)$, then at least one of the columns of $\gamma_n$ is eventually arbitrarily large (up to a subsequence), while the area of the parallelogram constructed out of these columns is $1$. As a consequence either one of these vectors diverge and the other remains bounded, or both vectors diverge and their directions get closer. This implies that for every diffuse faithful measure $\mu$ on $\mathbb{R}P^1$, the $\gamma_n$-translates of $\mu$ converge to a point measure, the Poisson map $C(\mathbb{R}P^1) \rightarrow l^\infty (\SL(2,\zz))$, $f \mapsto P_\mu (f) (\gamma)=\int f d (\gamma \mu)$ is an injective $*$-homomorphism and the right-action on its image is trivial; this shows that the action of $\SL(2,\zz) \times \SL(2,\zz)$ on $\partial_\beta \SL(2,\zz)$ is amenable. The same approach in the case of $\SL(3,\zz)$ immediately fails, since in this case the translates of a diffuse measure under a diverging sequence in $\SL(3,\zz)$ will not converge in general to a point mass, but measures supported on "great circles" are also allowed (see \cite{Fu} 4.6). The description of the Stone-{\v C}ech boundary of $\SL(3,\zz)$ we give in this work reflects this behaviour and leads to the study of two different parts of the boundary separatedly, depending on their "rank". 
 
 In Section \ref{s1} we establish some facts concerning Koopman representations of discrete groups. In Section \ref{s2} we give the description of the Stone-{\v C}ech boundary of $\SL(3,\mathbb{Z})$ as a fibered space and describe its decomposition according to the "rank". While the "rank $1$" part is an amenable space by construction, the "rank $2$" requires some considerations. The study of this region is performed in Section \ref{section4}. We consider preimages of some specific sets, which we call linear sets; these comprise sets of matrices with range containing a given line and kernel contained in a given plane, called (line, plane)-sets. We show that the quasi-invariant Radon measures annihilating these preimages have matrix coefficients which belong to a certain ideal. In the remaining part of Section \ref{section4} we show that the quasi-invariant Radon measures supported on these preimages are tempered. In order to show this we first prove amenability of the action on orbits of certain "dynamical atoms" in subsection \ref{secrat} and subsection \ref{secirr}; then we extend the temperedness result to (line, plane)-sets in subsection \ref{s5}.
 
 The construction given here can be generalized to other groups. Some explicit cases, such as $\SL(2,\mathbb{Z}[1/p])$, with $p$ prime, will be discussed in a subsequent paper (cfr. \cite{Ra} for a related topic).

  \section{Preliminaries}
  \label{s1}
  \subsection{Topological dynamics and Koopman representations}
  This section contains general results which will be used throughout the following chapters, in particular when discussing amenability of the left-right action of $\SL(3,\zz) \times \SL(3,\zz)$ on the rational and irrational parts of the boundary. 
   
  \begin{defn}
 Let $\Gamma$ be a countable discrete group acting on a set $X$. For every $Y \subset X$ we define the \textit{global stabilizer of $Y$} to be
 \begin{equation*}
 \Gamma_Y := \{ \gamma \in \Gamma \; | \; \gamma y \in Y \quad \forall y \in Y\}.
 \end{equation*}
 \end{defn}

 \begin{lem}
 \label{lemamencoset}
Let $\Gamma$ be a countable discrete group, $X$ a locally compact $\Gamma$-space, $Y \subset X$ a closed subset such that
\begin{itemize}
\item[(i)] $ \Gamma Y =X $,
\item[(ii)] for every $\gamma_1$, $\gamma_2 \in \Gamma$ we have
\begin{equation*}
\gamma_1 Y \cap \gamma_2 Y \neq \emptyset \qquad \Leftrightarrow \qquad \gamma_1 Y= \gamma_2Y.
\end{equation*}
\end{itemize}
Let $\Lambda $ be the global stabilizer of $Y$ (which is a group in this situation). If the action of $\Lambda$ on $Y$ is amenable, then the action of $\Gamma$ on $X$ is amenable.
\end{lem}
 \proof We shall use a convenient description of the dynamical system given by the action of $\Gamma$ on itself. Choose coset representatives $\{g_i\}$ for $\Gamma / \Lambda$. We can write $\Gamma = \Gamma / \Lambda \times \Lambda$; the left action of $\Gamma$ on itself becomes:
 \begin{equation*}
 \gamma (g_i \Lambda, \lambda) = (g_k \Lambda , \eta \lambda) \qquad \mbox{ for } \gamma g_i \Lambda = g_k \Lambda,
 \end{equation*}
 where $\eta \in \Lambda$ is uniquely determined by the requirement $\gamma g_i \Lambda = g_k \Lambda$.\\
 We can give a similar description of the action of $\Gamma$ on $X$. In order to see this, let $g_i Y \cap g_j Y \neq \emptyset$, then there is $y \in Y$ such that $g_i y \in g_j Y$ and so $g_j^{-1} g_i y \in Y$, which gives, by (ii), $g_j^{-1} g_i \in \Lambda$ and so $g_i \Lambda = g_j \Lambda$. Similarly, if $g_i Y \cap g_j Y = \emptyset$ then $g_i \Lambda \neq g_j \Lambda$. Hence the map from the set of translates $\{g_i Y\}$ to the coset space $\Gamma /\Lambda$ is well-defined and injective (it is surjective by construction). On the other hand, given $\gamma \in \Gamma$, we have $\gamma \Lambda = g_i \Lambda$ for some cose representative $g_i$; hence $g_i^{-1} \gamma \in \Lambda$, which gives (again by (ii)) $g_i Y = \gamma Y$. It follows from (i) that the translates of $Y$ are in bijection with the cosets $ g_i \Lambda $. We thus obtain the dynamical characterization of $X$ as $\Gamma / \Lambda \times Y$, where the action of $\Gamma$ is given by
 \begin{equation*}
 \gamma (g_i \Lambda , y) = (g_k \Lambda , \eta y) \qquad \mbox{ for } \gamma g_i \Lambda = g_k \Lambda.
 \end{equation*}
 Let now $\mu_\lambda : Y \rightarrow \mathcal{P}(\Lambda)$ be an approximate equivariant mean for the action of $\Lambda$ on $Y$ and consider $\eta_\lambda : \Gamma / \Lambda \times Y \rightarrow \mathcal{P} (\Gamma / \Lambda \times \Lambda)$, $\eta_\lambda ((g_i \Lambda, y))= \delta_{g_i \Lambda} \times \mu_\lambda (y)$. This is the desired approximate invariant mean. $\Box$
 

 Let $X$ be a locally compact space and $\Gamma$ a countable discrete group acting on it. For every quasi-invariant $\sigma$-finite measure $\mu$ on $X$ we can consider the associated Koopman representation $\pi_\mu$ of the group on $L^2 (X,\mu)$ given, for $\xi \in L^2 (X,\mu)$ and $\gamma \in \Gamma$, by
\begin{equation*}
(\gamma \xi)(x) = \sqrt{\frac{d \gamma \mu}{d \mu}} (x) \xi (\gamma^{-1} (x)),
\end{equation*}
 where $d\gamma \mu /d \mu$ is the Radon-Nikodym cocycle. In the same way, one can consider the Koopman representation $\pi_\mu^\rtimes$ of the full crossed product $C^*$-algebra $C_0 (X) \rtimes \Gamma$. Koopman representations have been studied in several places, see for example \cite{DuGr, BeGu} or \cite{Ke}, where the concept of temperedness for such representations is also discussed.\\
 By a Radon measure we mean a Borel measure which is outer regular (Borel sets are approximated in measure by open larger sets), inner regular (Borel sets are approximated in measure by smaller compact sets) and finite on compact sets. 

  \begin{lem}
  \label{wk}
  Let $\Gamma$ be a countable discrete group acting on a locally compact space $X$ and let $\mu$ be a quasi-invariant $\sigma$-finite Radon measure on $X$. If the action of $\Gamma$ on $X$ is amenable, the Koopman representation $C^* (\Gamma) \rightarrow \mathbb{B}(L^2 (X,\nu))$ factors through $C^*_{red}(\Gamma)$.
  \end{lem}

  \proof
Denoting by $\lambda$ the left regular representation of $\Gamma$, we have that $\lambda \otimes \pi_\mu \simeq \lambda \otimes 1$, where $1$ is the trivial representation on $L^2 (X,\mu)$. Hence we need to observe that the matrix coefficients for $\pi_\mu$ can be approximated by matrix coefficients for $\lambda \otimes \pi_\mu$. Let $\xi \in L^2 (X, \mu)$, since $\mu$ is Radon, we can approximate $\xi$ by compactly supported functions; then the matrix coefficients of $\xi \cdot f_i$ approximate the corresponding matrix coefficient of $\xi$, where $f_i$ is extracted from the net of \cite{An} Remark 4.10. $\Box$

 \begin{lem}
 \label{ideal}
 Let $C^*_\tau (\Gamma)$ be a $C^*$-completion of the group ring of a countable discrete group $\Gamma$. Let $X$ be a locally compact space on which $\Gamma$ acts and $\mu$ a quasi-invariant sigma-finite measure on $X$. Suppose that $\mu =\sum_{i=1}^\infty \mu_i$ with $\mu_i \perp \mu_j$ are quasi-invariant finite measures and that the Koopman representations $C^*( \Gamma) \rightarrow C^*_{\pi_{\mu_i} } (  \Gamma)$ factor through $C^*( \Gamma) \rightarrow  C^*_\tau (  \Gamma)$ for every $i \in \mathbb{N}$. Then the Koopman representation $C^*(  \Gamma) \rightarrow C^*_{\pi_\mu  }(  \Gamma)$ factors through $C^*_\tau (  \Gamma)$ as well.
 \end{lem}
 \proof
 Let $\phi$ be a state on $C^* (  \Gamma)$ associated to a function $\xi \in L^2_1 (\mu)$; for every $i$ denote by $\xi_i$ its component on $L^2 (\mu_i)$. For every $x \in C_c (\Gamma )$ we have $|( \pi_{\mu_i} (x) \xi_i, \xi_i)_{\mu_i}| \leq \| x\|_\tau \|\xi_i\|$. Hence $|(\pi_\mu (x) \xi, \xi)| =| \sum_i (\pi_{\mu_i} (x) \xi_i, \xi_i)_{\mu_i} |\leq \|x\|_\tau \sum_i \|\xi_i\|^2 = \|x\|_\tau \|\xi \|^2$. $\Box$\\

 \subsection{Exotic group $C^*$-algebras}
 
We recall the construction of exotic group $C^*$-algebras of a countable discrete group $\Gamma$ from $\Gamma \times \Gamma$-invariant ideals in $l^\infty (\Gamma)$, as developed in \cite{BrGu}. Fix such a group $\Gamma$ and an ideal $I\subset l^\infty (\Gamma)$ containing $c_c (\Gamma)$. 
\begin{defn}[\cite{BrGu} Definition 2.2]
\label{def-irep}
An {\it $I$-representation} of $C^*(\Gamma)$ is a unitary representation on a Hilbert space $H$ admitting a dense set of vectors $K \subset H$ such that for every $\xi, \eta \in K$ the associated matrix coefficient $\{\gamma \mapsto \braket{ \gamma \xi, \eta}\}$ is an element of $I$.
\end{defn}
 Define a $C^*$-norm on the group ring by $\|x\|_I = \sup\{ \|\pi (x)\| \; | \; \pi \mbox{ is an $I$-representation }\}$; the completion of $\cc [\Gamma]$ with respect to $\| \cdot \|_I$ is the exotic group $C^*$-algebra associated to the ideal $I$ and is denoted by $C^*_I (\Gamma)$. A similar construction can be developed for the case of crossed products. For more information concerning exotic group $C^*$-algebras and exotic crossed products see for example \cite{RuWi,BuEcWi}. The following will be useful in reducing our study to a particular region of the Stone-{\v C}ech boundary of $\SL(3,\mathbb{Z})$.

\begin{lem}
\label{inv-ideal}
Let $\Gamma$ be a discrete countable group, $\alpha : \Gamma \rightarrow \Gamma$ a group automorphism and $I \subset l^\infty (\Gamma)$ an $\alpha$-invariant ideal, i.e. both $f \circ \alpha$ and $f \circ \alpha^{-1}$ belong to $I$ for every $f \in I$. Then for every $x \in C^* (\Gamma)$ we have
\begin{equation*}
\| x\|_I = \| \alpha (x) \|_I,
\end{equation*}
where we consider $\alpha$ as a $*$-automorphism of $C^* (\Gamma)$. In particular, the kernel of the quotient map $C^*(\Gamma) \rightarrow C^*_I (\Gamma)$ is $\alpha$-invariant.
\end{lem}
\proof Let $\pi$ be an $I$-representation and $\xi$ a vector in the corresponding Hilbert space whose associated matrix coefficient $\phi_\xi$ belongs to $I$. Since $I$ is $\alpha$-invariant, also the matrix coefficient $\phi_\xi \circ \alpha$ belongs to $I$ and so $\pi \circ \alpha$ is a $c_0$-representation. The same is true for $\pi \circ \alpha^{-1}$; hence every $I$-representation is of the form $\pi \circ \alpha$, where $\pi$ is an $I$-representation and for every $x \in C^*(\Gamma)$ we have
\begin{equation*}
\| x \|_I = \sup  \{\| \pi \circ \alpha (x) \|  \; | \; \pi \mbox{ is an $I$-representation } \}= \| \alpha (x) \|_I,
\end{equation*}
which is the desired equality. $\Box$\\

\begin{defn}
\label{def-c0}Let $\Gamma$ be a discrete countable group. We define the translation invariant ideal in $l^\infty (\Gamma \times \Gamma)$ given by $c'_0 (\Gamma \times \Gamma):= \{ f \in l^\infty (\Gamma \times \Gamma) \; | \; f(x_n, y_n) \rightarrow 0\; \mbox{ whenever both } x_n\rightarrow \infty \mbox{ and } y_n \rightarrow \infty\}$.\end{defn}


\subsection{Proximal pieces}
The concept of proximal (boundary) piece for a discrete countable group was introduced in \cite{BIP} as a generalization of the notion of boundary which is small at infinity in order to study rigidity properties for groups which do not satisfy any negative curvature condition (see the Introduction in \cite{BIP}).

\begin{defn}[\cite{BIP} Definition 3.7]
Let $\Gamma$ be a discrete countable group acting on a compact space $K$ and $\eta$ a probability measure on $K$, which is not invariant under $\Gamma$. The \textit{proximal (boundary) piece} associated to $(K, \eta)$ is the subset of $\partial_\beta \Gamma$ of elements $\omega$ satisfying $\lim_{g \rightarrow \omega} ( (gh)\cdot \eta - g \cdot \eta) =0$ in the $w^*$-topology, for every $h \in \Gamma$. It is denoted by $\partial_\eta \Gamma$.
\end{defn}

Proximal boundary pieces are closed $\Gamma \times \Gamma$-invariant subsets of $\partial_\beta \Gamma$, the condition that the measure $\eta$ is not invariant guarantees that that proximal pieces are actually subsets of the boundary $\partial_\beta \Gamma$, rather than the whole Stone-\v{C}ech compactification $\Delta_\beta \Gamma$. In the case $\Gamma = \SL(d, \mathbb{Z})$, $d \geq 2$, natural examples of proximal pieces arise from actions on Grassmanians, endowed with the unique $SO (d)$-invariant probability measure; a proximal piece of this sort admits a concrete description in terms of the multi-index defining the Grassmanian (\cite{BIP} Chapter 6). We will eventually be interested in the proximal piece for $\SL(3,\mathbb{Z})$ associated to the full Grassmanian corresponding to the multi-index $(1,2,3)$, or equivalently, to the Borel subgroup $B(\mathbb{R}) \subset \SL(3,\mathbb{R})$ of upper triangular matrices. Note that in this case the pair $(G(\mathbb{R}/B(\mathbb{R}), \eta)$ is the Poisson boundary of $(\SL(d,\mathbb{R}), \mu)$, with $\mu$ a spherical measure (\cite{Fu} Theorem 4.1). We will denote this particular proximal piece by $\partial \SL(3,\mathbb{Z})_0$.

\subsection{Calkin states for discrete countable groups}

Let $\Gamma$ be a discrete countable group, $X$ be a compact Hausdorff space and $\mu$ a $\Gamma$-quasi-invariant Radon measure on $X$. For every $\gamma \in \Gamma$ we consider the representation $\pi_{\gamma, \mu}$ of $C(X)$ on $L^2(\mu)$ given by $(\pi_{\gamma, \mu}(f) \xi)(x)= f(\gamma^{-1} (x)) \xi (x)$.

\begin{lem}
\label{lem2.4.1}
Let $\Gamma$ be a discrete group, $X$ be a compact Hausdorff space and $\mu$ a $\Gamma$-quasi-invariant Radon measure on $X$. The Koopman representation is the only group homomorphism $u : \Gamma \rightarrow U(L^2 (\mu))$ satisfying:
\begin{itemize}
\item for every $\gamma \in \Gamma$, $u_\gamma$ is a unitary $C(X)$-intertwining between the representations $(L^2(\mu),\pi_e)$ and $(L^2(\mu), \pi_\gamma)$;
\item for every $\gamma \in \Gamma$, $u_\gamma$ sends positive elements in $L^2(\mu)$ to positive elements.
\end{itemize}
\end{lem}
\proof Note that $C(X) \cap \mathbb{B}(L^2 (\mu))=L^\infty (\mu)$, hence if $U$ and $V$ are two unitaries in $\mathbb{B}(L^2(\mu))$ which implement the same automorphism of $C(X)$, they differ by a unitary in $L^\infty (\mu)$. The positivity condition then entails the result. $\Box$

\begin{lem}
\label{lem2.4.2}
Let $X$ be a compact Hausdorff space, $\mu$ a Radon measure on $X$ and $V \subset C(X)$ a linear space which is closed under passage to the absolute value and whose image in $L^2 (\mu)$ is dense. Then every positive function in $L^2 (\mu)$ is approximated by positive elements coming from $V$.
\end{lem}
\proof 
Let $\xi \in L^2(\mu)_+$ and $\epsilon >0$. There is $a \in V$ such that $\|\xi-a\|^2_2 < \epsilon$. Let $F:= \{ x \in X \; | \; a(x) >0\}$. Then $\int_X \xi^2 + \int_{X \backslash F} |\xi-a|^2 \leq \int_F |\xi - a|^2 + \int_{X \backslash F} |\xi - a|^2 < \epsilon$. Define $b=(a + |a|)/2$. Then $b \geq 0$ and $\int_X |\xi - b|^2 = \int_F \xi^2 - \int_{X \backslash F} |\xi - a|^2 < \epsilon$. $\Box$

If $A$ is a unital $C^*$-algebra and $\omega \in \partial_\beta \mathbb{N}$, we denote by $A_\omega := l^\infty(\mathbb{N}, A)/I_\omega$, where $I_\omega := \{ a \in l^\infty (\mathbb{N},A) \; | \; \lim_{n\rightarrow \omega} \|a\| =0\}$, the ultraproduct associated to $\omega$. If $\Gamma$ is a countable discrete group, for every $\omega \in \partial_\beta \mathbb{N}$, there is a surjection $\psi_\omega : \sigma(l^\infty(\Gamma)_\omega) \rightarrow \Delta_\beta \Gamma$ equivalent to the diagonal embedding $l^\infty(\Gamma) \rightarrow l^\infty (\Gamma)_\omega$, $f \rightarrow (f_n)=(f)$.

\begin{thm}
\label{thm2.4.3}
Let $\Gamma$ be a discrete countable group and denote by $\pi_{\mathcal{C}}$ the left-right representation of $C^* (\Gamma \times \Gamma)$ in the Calkin algebra $\mathbb{B}(l^2 (\Gamma)) / \mathbb{K} (l^2 (\Gamma))$. Let $C^*_\tau (\Gamma \times \Gamma)$ be a quotient of $C^*(\Gamma \times \Gamma)$ such that for every $\omega \in \partial_\beta \mathbb{N}$ the Koopman representation associated to every $\Gamma \times \Gamma$-quasi-invariant Radon measure on $\psi_\omega^{-1} (\partial_\beta \Gamma) $ is a $C^*_\tau (\Gamma \times \Gamma)$-representation. Then $\pi_{\mathcal{C}} (C^*(\Gamma \times \Gamma))$ factors through $C^*_\tau (\Gamma \times \Gamma)$.
\end{thm}
\proof In virtue of \cite{Wi} 2.3.5 given a dense sequence $\{\xi_n\}_{n \in \mathbb{N}}$ of norm-one vectors in the unit ball of $l^2(\Gamma)$, every state $\phi$ on $\mathbb{B}(l^2 (\Gamma))/\mathbb{K}(l^2(\Gamma))$ is of the form $\phi (x) = \lim_{n \rightarrow \omega} \braket{x\xi_n,\xi_n}$ for some $\omega \in \partial_\beta \mathbb{N}$. Note that the fact that $\phi$ defines a state on the Calkin algebra implies that the $\xi_n$'s converge weakly to zero with respect to $\omega$. So let be given such a sequence $\{\xi_n\}$ and such an ultrafilter $\omega \in \partial_\beta \mathbb{N}$, corresponding to a state $\phi$. We want to show that the restriction of $\phi$ to $\pi_{\mathcal{C}}(C^*(\Gamma \times \Gamma))$ is the $w^*$-limit of states associated to Koopman representations of $\Gamma \times \Gamma$ on $\sigma (l^\infty(\Gamma)_\omega)$ given by $\Gamma \times \Gamma$-quasi-invariant finite Radon measures, supported on $\psi_\omega^{-1} (\partial_\beta \Gamma)$. For it is enough to show that there are a sequence of finite quasi-invariant Radon measures $\mu_k$ on $\psi_\omega^{-1} (\partial_\beta \Gamma)$ and a uniformly bounded sequence (in the respective $\| \cdot \|_2$-norms) of vectors $\eta_k \in L^2 (\mu_k)$ such that $\phi (\gamma) = \lim_{k \rightarrow \infty} \braket{\gamma \eta_k, \eta_k}_{L^2 (\mu_k)}$ for every $\gamma \in \Gamma \times \Gamma$. We will show the stronger fact that there are a unique $\Gamma \times \Gamma$-quasi-invariant measure $\tilde{\mu}$ and a vector $\xi \in L^2(\tilde{\mu})$ such that for every $\gamma \in \Gamma \times \Gamma$ we have $\braket{ \gamma \xi, \xi}_{L^2(\tilde{\mu})}= \lim_{n \rightarrow \omega} \braket{\gamma \xi_n, \xi_n}$.\\
Let $\{\gamma_i\}_{i=0}^\infty$ be an enumeration of $\Gamma \times \Gamma$ with $\gamma_0 =e$. For every $n \in \mathbb{N}$ let $\tilde{\xi}_n (x) := \sum_{i=0}^\infty \frac{1}{2^i} |\xi_n| (\gamma_i^{-1} x)$. The $\tilde{\xi}_n$'s constitute a uniformly bounded sequence in $l^2(\Gamma)$, which converges weakly to zero with respect to $\omega$. Let $\tilde{\mu}$ be the Radon measure on $\sigma(l^\infty(\Gamma)_\omega)$ induced by the $\tilde{\xi}_n$'s by $\tilde{\mu} (f_n) = \lim_{n\rightarrow \omega} \braket{f_n \tilde{\xi}_n ,\tilde{\xi}_n}$, where $(f_n) \in l^\infty(\Gamma)_\omega$. Note that this measure is supported on $\psi_\omega^{-1} (\partial_\beta \Gamma)$. First of all we shall see that $\tilde{\mu}$ is $\Gamma \times \Gamma$-quasi-invariant. For it is enough to show that for every $\gamma \in \Gamma \times \Gamma$ and every $\epsilon >0$ there is $\delta >0$ such that for every $0 \leq f \leq 1$ in $l^\infty(\Gamma)_\omega$ satisfying $\tilde{\mu}(f) <\delta$ we have $\tilde{\mu}(f \circ \gamma) < \epsilon$. This is in turn equivalent to the claim that if $0 \leq f_k \leq 1$ is a sequence in $l^\infty (\Gamma)_\omega$ such that $\lim_k \tilde{\mu}(f_k) =0$, then $\lim_k \tilde{\mu} (f_k \circ \gamma)=0$ for every $\gamma \in \Gamma \times \Gamma$.\\
Now let $f_k$ be such a sequence, then 
\begin{equation*}
\begin{split}
0&=\lim_k \lim_{n \rightarrow \omega} \sum_x (f_k)_n (x) \tilde{\xi}_n(x)^2 = \lim_k \lim_{n \rightarrow \omega}\sum_x \sum_{i,j} 2^{-i-j} (f_k)_n(x) |\xi_n| (\gamma_i^{-1} x) |\xi_n|(\gamma_j^{-1} (x)) \\
&= \lim_k \lim_{n \rightarrow \omega} \sum_{i,j} 2^{-i-j} \braket{(f_k)_n |\xi_n|\circ \gamma_i^{-1}, |\xi_n| \circ \gamma_j^{-1}},
\end{split}
\end{equation*}
from which it follows that for every $i,j$ we have $\lim_k \lim_{n\rightarrow \omega} \braket{ (f_k)_n |\xi_n| \circ \gamma_i^{-1}, |\xi_n| \circ \gamma_j^{-1}} =0$. Hence, given $\gamma \in \Gamma \times \Gamma$, for every $i,j$ we have $\lim_k \lim_{n\rightarrow \omega} \braket{ (f_k)_n |\xi_n| \circ( \gamma_i^{-1} \gamma^{-1}), |\xi_n| \circ (\gamma_j^{-1} \gamma^{-1})} =0$, which gives 
\begin{equation*}
\lim_k \tilde{\mu} ((f_k)_n \circ \gamma)=\lim_k \lim_{n \rightarrow \omega} \sum_{i,j} 2^{-i-j} \braket{(f_k)_n |\xi_n|\circ (\gamma_i^{-1} \gamma^{-1}), |\xi_n| \circ (\gamma_j^{-1} \gamma^{-1})}=0.
\end{equation*}
Let $W$ be the linear subspace of $\prod_n l^\infty (\Gamma)$ defined as 
\begin{equation*}
\begin{split}
W:=\{&f=(f_n) \in \prod_n l^\infty (\Gamma) \; | \; \mbox{for every $\gamma \in \Gamma \times \Gamma$ there is $M_\gamma >0$ satisfying } \\
&|f_n\circ \gamma^{-1}| \leq M_\gamma \tilde{\xi}_n \mbox{ for every $n$ }\}.
\end{split}
\end{equation*}
Consider classes in $L^2(\tilde{\mu})$ of elements in $\prod_n l^\infty (\Gamma)$ of the form $(f_n/\tilde{\xi}_n)$ with $f=(f_n) \in W$. Denote the resulting linear subspace of $L^2(\tilde{\mu})$ by $W_{\tilde{\mu}}$. This is a dense linear subspace of $L^2(\tilde{\mu})$, which is also an $l^\infty(\Gamma)_\omega$-module, in fact, if $f=(f_n) \in L^2(\tilde{\mu})$ is represented by an element of $l^\infty(\Gamma)_\omega$, then for every $\epsilon >0$ there is $N=N_\epsilon$ such that $\sum_{i=N+1}^\infty 2^{-i} < \sqrt{\epsilon} \|f\|^{-1}$ and so
\begin{equation*}
\begin{split}
\sum_x|f_n(x) &- \sum_{i=0}^N f_n(x) \frac{|\xi_n|(\gamma_i^{-1}(x))}{\tilde{\xi}_n(x)}|^2 \tilde{\xi}_n(x)^2 = \sum_x|f_n(x)[\tilde{\xi}_n (x) - \sum_{i=0}^N  |\xi_n|(\gamma_i^{-1}(x))]|^2\\
&=\sum_x|f_n(x)|^2 |\sum_{i=N+1}^\infty 2^{-i} |\xi_n|(\gamma_i^{-1}x)|^2 \leq \|f_n\|^2 \sum_{i,j=N+1}^\infty 2^{-i-j} \sum_x |\xi_n|(\gamma_i^{-1}(x))|\xi_n| (\gamma_j^{-1} (x))\\ &\leq \|f_n\|^2 \sum_{i,j=N+1}^\infty 2^{-i-j} 
\end{split}
\end{equation*}
and so
\begin{equation*}
\|f - f\sum_{i=0}^N |\xi| \circ \gamma_i^{-1}/\tilde{\xi}\|_2^2 < \epsilon.
\end{equation*}
 For every $\gamma \in \Gamma \times \Gamma$ define the map $u_\gamma$ on $W\subset \prod_n l^\infty (\Gamma)$ given by $f/\tilde{\xi}_n \mapsto f\circ \gamma^{-1} /\tilde{\xi}_n$. We want to check that this defines a linear map on $W_{\tilde{\mu}}\subset L^2(\tilde{\mu})$. For suppose that $f$, $g \in W \subset \prod_n l^\infty (\Gamma)$ are such that $(f_n/\tilde{\xi}_n)=(g_n /\tilde{\xi}_n)$ in $L^2 (\tilde{\mu})$. Then $\lim_{n \rightarrow \omega} \sum_x |f_n(x)-g_n(x)|^2 =0$ and so $\lim_{n \rightarrow \omega} \sum_x |f\circ \gamma^{-1}(x) - g\circ \gamma^{-1} (x)|^2 = \|f_n\circ \gamma^{-1}/\tilde{\xi}_n - g_n \circ \gamma^{-1} /\tilde{\xi}_n\|^2 =0$. Hence $u_\gamma$ is well defined on $W_{\tilde{\mu}}$ and extends to a unitary intertwining between $\pi_e$ and $\pi_\gamma$. The dense subspace $W_{\tilde{\mu}}$ is closed under passage to the absolute value and so it follows from Lemma \ref{lem2.4.1} and Lemma \ref{lem2.4.2} that the extension of $u_\gamma$ to $L^2(\tilde{\mu})$ is the Koopman operator. Hence, taking $\xi=(\xi_n)$, we can compute
\begin{equation*}
\begin{split}
\braket{\gamma (\xi/\tilde{\xi}), \xi/\tilde{\xi}}_{\tilde{\mu}}&= \lim_{n \rightarrow \omega} \braket{\xi_n \circ \gamma^{-1}/\tilde{\xi}_n, \xi_n /\tilde{\xi}_n} = \lim_{n \rightarrow \omega} \sum_x \frac{\xi_n \circ \gamma^{-1}(x) \xi_n(x) }{\tilde{\xi}_n(x)^2} \tilde{\xi}_n (x)^2\\ &= \lim_{n\rightarrow \omega} \braket{\xi_n \circ \gamma^{-1}, \xi_n}.
\end{split}
\end{equation*}
The result follows. $\Box$

\section{A suitable boundary}
\label{s2}
From now on we focus on the case of $\SL(3,\zz)$. This section contains a description of the Stone-{\v C}ech boundary of this group as a fibered space. This will be used in the following chapters, where we study specific preimages in more details. 

Let $G=\SL(3 )$ and $Y$ be the compact space obtained as the cartesian product of the space of projective $3 \times 3$ matrices $\mathbb{P}M_3 =(M_3 \backslash \{0\}) /\mathbb{R}_{\neq 0}$ with itself, $Y=\mathbb{P}M_3  \times \mathbb{P}M_3 $. This space is endowed with an action of $G(\rr) \times G(\rr) $ given, for $(g_1, g_2) \in G(\rr) \times G(\rr)$, by $(g_1,g_2) (m_1,m_2) =(g_1 m_1  g_2^{-1} / \mathbb{R}_{\neq 0}, g_2 m_2 g_1^{-1}/\mathbb{R}_{\neq 0 })$. The $G(\mathbb{Z}) \times G(\mathbb{Z})$- equivariant map $\phi : G(\mathbb{Z})\rightarrow Y$ given by $\gamma \mapsto (\gamma /\mathbb{R}_{\neq 0}, \gamma^{-1}/\mathbb{R}_{\neq 0})$ is continuous and extends to a continuous $G(\mathbb{Z}) \times G(\mathbb{Z})$-equivariant map from the Stone-{\v C}ech compactification $\Delta_\beta G(\mathbb{Z})$ of $G(\mathbb{Z})$ to $Y$. Let $Y_1 \subset Y$ be the subset of pairs of projective rank $1$ matrices and $Y_2 \subset Y$ be the subset of pairs of projective matrices $(m_1,m_2)$ such that either $\rank (m_1)=1$, $\rank(m_2)=2$ or $\rank (m_1)=2$, $\rank (m_2)=1$ and $m_1 m_2 =m_2 m_1 =0$ (matrix multiplication).

 To every action of a countable discrete group on a compact space and probability measure on it, it is possible to associate a proximal boundary piece (\cite{BIP}, 3).

\begin{prop}
\label{prop1}
Let $\partial G(\mathbb{Z})_0$ be the proximal boundary piece associated to the action of $G(\mathbb{Z})$ on the Poisson boundary $(G(\mathbb{R})/B(\mathbb{R}), \eta)$ of $G(\mathbb{R})$. Then  
\begin{equation*}
\phi^{-1} (Y_1)= \{ \omega \in \partial_\beta G(\mathbb{Z}) \; | \; \phi (\omega) \in Y_1\} \subseteq \partial G(\mathbb{Z})_0 .
\end{equation*}
Moreover, $\phi (\partial_\beta G(\mathbb{Z})) \subset Y_1 \cup Y_2$. $\phi^{-1}(Y_1)$ is a closed subset of $\partial_\beta G(\mathbb{Z})$ and $\phi^{-1}(Y_2)$ is open. The action of $G(\zz) \times G(\zz)$ on $ Y_1 $ is amenable.
\end{prop}
 \proof Let $g_n$ be a diverging sequence in $G(\mathbb{Z})$, then $\det(g_n /\| g_n\|) = \|g_n\|^{-3}$ tends to zero, as well as $\det(g_n^{-1}/\|g_n^{-1}\|)$. Hence $\phi|_{\partial G(\mathbb{Z})}$ takes values in pairs of non-invertible matrices. Moreover, since the product of matrices is continuous on $M_3(\mathbb{R})$, it follows that for $\lim_n (g_n/\|g_n\|, g_n^{-1}/\|g_n^{-1}\|)=(m_1,m_2 )$, we have $m_1 m_2 = \lim_n g_n g_n^{-1} /(\|g_n\| \|g_n^{-1}\|) =0$ and so $\phi (\partial G(\mathbb{Z})) \subset Y_1 \cup Y_2$. Let now $\omega \in \partial_\beta G(\mathbb{Z})$ be such that $\phi (\omega)$ belongs to $Y_1$. For every $g \in G(\mathbb{Z})$ let $s_{i }(g)$, $i=1,2,3$ be the eigenvalues of $\sqrt{g^t g }$, ordered in such a way that $s_{1 }(g) \geq s_{2 }(g) \geq s_{3 }(g)$ (cfr. \cite{BIP}, 6.1). Consider the functions $  f_{1,2}: g \mapsto  s_1(g)/s_2(g) $ and $ f_{2,3} : g \mapsto  s_2(g) / s_3(g) $. Suppose that $f_{1,2} (\omega) = c < \infty$. Then there is a sequence $g_n$ with $\lim_n f_{1,2} (g_n) =c$ and $\lim_n \phi (g_n) \in Y_1$. Following \cite{BIP} Lemma 6.1 for every $n \in \mathbb{N}$ we can write 
 \begin{equation*}
g_n=a_n s_n b_n = a_n\left( \begin{array}{ccc}	s_{1,n}	&	0	&	0	\\
							0		&	s_{2,n}	&	0	\\
								0	&	0		&	s_{3,n}\end{array}\right)b_n,
\end{equation*}
where $a_n, b_n \in \SO(3,\mathbb{R})$ and $s_{i,n} = s_i (g_n)$ for $i \in \{1,2,3\}$. Up to taking a subsequence, we can suppose that $a_n \ra a$, $b_n \ra b$ in $\SO(3,\mathbb{R})$, thus also $\lim_n \phi (s_n)$ exists and $\lim_n \phi (a_n s_n b_n)= \lim_n a_n \phi (s_n) b_n= a\lim_n \phi(s_n) b$; but the condition $\lim_n f_{1,2} (g_n) =c$ entails $\rank (\phi (\omega))=\rank(\lim_n s_n/\|s_n\|))=2$. In a similar way, if $f_{2,3} (\omega) =c < \infty$, then there is a sequence $g_n$ with $\lim_n \phi(g_n) = \phi(\omega)$ and $\lim_n f_{2,3}(g_n)=f_{2,3}(\omega) =c< \infty$, but then $g_n^{-1}/\|g_n^{-1}\|$ converges to a rank-$2$ matrix, which is impossible. Hence $\phi^{-1} (Y_1) \subseteq \partial G(\mathbb{Z})_0$. Note now that the set of pairs of norm-$1$ matrices $(m_1, m_2)$ in $M_3(\mathbb{R}) \times M_3 (\mathbb{R})$ satisfying $m_1  m_2 = m_2 m_1 =0$ is closed in the product topology of $(M_3 \backslash \{0\}) \times (M_3 \backslash \{0\})$, being defined by the zeros of two polynomials. The quotient map $(M_3 \backslash \{0\}) \times (M_3 \backslash \{0\}) \rightarrow \mathbb{R}(M_3) \times \mathbb{P}(M_3)$ is open, being obtained from an action of $\mathbb{R}^\times \times \mathbb{R}^\times$; hence $Y_1$ is closed. The amenability of the action of $G(\zz) \times G(\zz)$ on $Y_1$ can be shown in the following way: consider the map $\psi: Y_1 \rightarrow G(\rr)/B(\rr) \times G(\rr) / B(\rr)$ given by $\psi (m_1,m_2) = ( \ran (m_1) \subset \ker (m_2), \ran (m_2) \subset \ker (m_1))$, where we see the inclusions $\ran(m) \subset \ker (m')$ as flags; this map is $G(\rr) \times G(\rr)$-equivariant. $\Box$
 
 Note that the subsets of $Y_2$ the form $M_{3,2} = \{ m=(m_1, m_2) \in Y_2 \; | \; \rank (m_1)=2,\; \rank (m_2)=1\}$ and $M_{3,2}' = \{ m=(m_1, m_2) \in Y_2 \; | \; \rank (m_1)=1,\; \rank (m_2)=2\}$ are $G(\mathbb{Z}) \times G(\mathbb{Z})$-invariant.  
 \begin{lem}
 \label{lem2}
 The spaces $M_{3,2}$, $M_{3,2}'$ are homogeneous spaces for the action of $G(\rr) \times G(\rr)$, with $M_{3,2}$ $G(\rr) \times G(\rr)$-homeomorphic to $(G(\rr) \times G(\rr))/P(\rr)$, where 
 \begin{equation*}
 P(\mathbb{R})=\{ (g_1, g_2) \in G(\rr) \times G(\rr) \; | \; \exists g \in \GL(2,\mathbb{R}), t_1,t_2,s_1,s_2 \in \mathbb{R} ,  \alpha \in \mathbb{R}_{> 0}\; : 
 \end{equation*}
 \begin{equation*}
 g_1=\left(\begin{array}{ccc} g	&	&	t_1\\
 		&	&	t_2\\
		0	&	0	&	\det(g)^{-1}	\end{array}\right), \quad g_2=\left(\begin{array}{ccc} \alpha g 	&	&	0\\
			&	&	0\\
			s_1	&	s_2	&	\det( \alpha g)^{-1}	\end{array}\right) \} 
			\end{equation*}
and $M_{3,2}' \simeq (G(\rr) \times G(\rr))/P'(\rr)$, where
 \begin{equation*}
  P'(\mathbb{R})=\{ (g_1, g_2) \in G(\rr) \times G(\rr) \; | \; \exists g \in \GL(2,\mathbb{R}), t_1,t_2,s_1,s_2 \in \mathbb{R}, \alpha \in \mathbb{R}_{> 0} \; : 
 \end{equation*}
  \begin{equation*}
  g_1=\left(\begin{array}{ccc} g	&	&	0\\
 		&	&	0\\
		t_1 	&	t_2	&	\det(g)^{-1}	\end{array}\right), \quad g_2=\left(\begin{array}{ccc} \alpha g 	&	&	s_1\\
			&	&	s_2\\
			0	&	0	&	\det( \alpha g)^{-1}	\end{array}\right) \} 
 			\end{equation*}
  \end{lem}
 \proof Let $x \in M_{3,2}$ be the element
 \begin{equation*}
x= \left(\left(\begin{array}{ccc}	1	&	0	&	0\\	
 					0	&	1	&	0\\
					0	&	0	&	0	\end{array}\right),\left(
					\begin{array}{ccc}	0	&	0	&	0\\
									0	&	0	&	0\\
									0	&	0	&	1\end{array}\right)\right).
\end{equation*}
The stabilizer of this point is $P(\mathbb{R})$. Since the second component of the elements in $M_{3,2}$ is uniquely determined by the first one, we need to show that the action of $G(\mathbb{R}) \times G(\mathbb{R})$ on th space of $3\times 3$-projective matrices of rank $2$ is transitive. Let $m$ be such a matrix. Using the singular value decomposition we can write $m= \pm k_1 m_d k_2$, where $k_1$, $k_2 \in \SO(3,\mathbb{R})$ and $m_d$ is a diagonal matrix of the form
\begin{equation*}
m_d=\left( \begin{array}{ccc}	\lambda_1 	&	0	&	0\\
						0	&	\lambda_2		&	0\\
						0	&	0	&	0	\end{array}\right),
\end{equation*}
with $\lambda_1 \geq \lambda_2 >0$. Then $m_d$ is clearly in the orbit of $x$ and the action is transitive. A similar argument shows that $M'_{3,2} = (G(\mathbb{R}) \times G(\mathbb{R}))/P' (\mathbb{R})$. $\Box$\\   

The spaces $M_{3,2}$ and $M'_{3,2}$ are topologically conjugated when we consider on $M_{3,2}$ the natural action of $G(\mathbb{R}) \times G(\mathbb{R})$ and on $M'_{3,2}$ the composition of the natural action with the group automorphism $ \alpha: (g_1,g_2) \mapsto (g_2,g_1)$. We denote this topological conjugacy by $\bar{\alpha}$. At the level of the fibers, we have the following

\begin{prop}
\label{prop-inverse}
Let $\inv : \Delta_\beta G(\mathbb{Z}) \rightarrow \Delta_\beta G(\mathbb{Z})$ be the extension of the inverse map to the Stone-{\v C}ech compactification of $G(\mathbb{Z})$. Then $\inv$ sets up a topological conjugacy between $\phi^{-1}(M_{3,2})$ endowed with the natural action of $G(\mathbb{Z}) \times G(\mathbb{Z})$ and $\phi^{-1} (M'_{3,2})$ endowed with the $\alpha$-twisted action of $G(\mathbb{Z}) \times G(\mathbb{Z})$. The result remains true if we equip $\phi^{-1} ( M_{3,2})$ with the $\alpha$-twisted action and $\phi^{-1}(M'_{3,2})$ with the natural action.
\end{prop} 
\proof It is clear that the map $\inv$ interchanges the natural action with the $\alpha$-twisted action. Let $\omega \in \phi^{-1} (M_{3,2})$; we want to check that $\phi (\inv (\omega))$ belongs to $M'_{3,2}$. For let $\phi (\omega)=(m_2 , m_1)$ and let $U  \times V$ be a product of open sets in $Y$. The set $\{ \gamma \in G(\mathbb{Z}) \;  | \; \gamma \in U, \; \gamma^{-1} \in V \}$ belongs to $\omega$. Hence $\{ \gamma \in G(\mathbb{Z}) \; | \; \gamma \in V, \; \gamma^{-1} \in U\}$ belongs to $\inv (\omega)$. Since $U$ and $V$ are arbitrary, we obtain $\inv (\omega) = (m_2, m_1)$. The result follows. $\Box$

To summarize, we have the following commutative diagram of $G(\mathbb{Z}) \times G(\mathbb{Z})$-equivariant maps (we denote by $\cdot$ the natrual actions and by $\cdot_\alpha$ the $\alpha$-twisted actions)
\begin{equation*}
\begin{array}{ccc}	(\phi^{-1}(M_{3,2}), \cdot) &\overset{\inv}{\simeq} & (\phi^{-1} (M'_{3,2}), \cdot_{\alpha})\\
				  \downarrow 	&	&	\downarrow\\
				(M_{3,2}, \cdot) & \overset{\bar{\alpha}}{\simeq} & (M'_{3,2}, \cdot_{\alpha})
				\end{array}.
				\end{equation*}

For our purposes, we want to verify that a Koopman representation associated to a measure on $\phi^{-1}(M_{3,2})$ is a representation of $C^*_{red}(G(\mathbb{Z}) \times G(\mathbb{Z}))$ if and only if the same is true for the corresponding Koopman representation on $\phi^{-1}(M'_{3,2})$ (under $\bar{\alpha}$), endowed now with the {\it natural} action of $G(\mathbb{Z}) \times G(\mathbb{Z})$. In virtue of Proposition \ref{prop-inverse} this is in turn equivalent to the fact that the Koopman representation on $\phi^{-1}(M_{3,2})$ with the $\alpha$-twisted action is a quotient of $C^*_{red}(G(\mathbb{Z}) \times G(\mathbb{Z})) $. 

\begin{prop}
\label{pinv-koop}
Let $\mu$ be a $\sigma$-finite measure on $\phi^{-1}(M_{3,2})$. Consider the Koopman representation $\pi_\mu$ associated to the natural action and the Koopman representation $\pi_{\mu}^\alpha$ associated to the $\alpha$-twisted action. Then $\pi_\mu$ factors through $C^*_{red}(G(\mathbb{Z}) \times G(\mathbb{Z})) $ if and only if $\pi_\mu^\alpha$ does. Moreover, $\pi_\mu$ is a $c'_0 (G(\mathbb{Z}) \times G(\mathbb{Z}))$-representation if and only if $\pi_\mu^\alpha$ is.
\end{prop}
\proof For every $(g_1, g_2) \in G(\mathbb{Z}) \times G(\mathbb{Z})$ and $\xi \in L^2 (\mu)$ we have $\pi_\mu^\alpha (g_1,g_2) \xi = \alpha (g_1,g_2) \xi = (g_2,g_1) \xi$. Hence $\pi_\mu^\alpha = \pi_\mu \circ \alpha$ and $\ker (\pi_\mu^\alpha) = \alpha^{-1} (\ker (\pi_\mu))$. Let $K$ be the kernel of the quotient map $C^*(G(\mathbb{Z}) \times G(\mathbb{Z}))\rightarrow C^*_{red}(G(\mathbb{Z}) \times G(\mathbb{Z})) $. Since $c_c (G(\mathbb{Z}) \times G(\mathbb{Z}))$ is $\alpha$-invariant, it follows from Lemma \ref{inv-ideal} that also $K$ is $\alpha$-invariant. Hence, $K \subset \ker (\pi_\mu) \cap \ker (\pi_\mu^\alpha)$. 
Suppose now that $H_0$ is a dense set of vectors in $L^2(\mu)$ whose associated matrix coefficients for $\pi_\mu$ belong to $c'_0 (G(\mathbb{Z}) \times G(\mathbb{Z}))$, then the same is true for $\pi_\mu^\alpha$ in place of $\pi_\mu$. $\Box$

 

 \section{Measures and linear sets in $M_{3,2}$}
 \label{section4}
 
We want to study the topology of the image of $C^*(G(\zz) \times G(\zz))$ in the Calkin algebra of $G(\zz)$. In virtue of Theorem \ref{thm2.4.3} we are led to consider the Koopman representations associated to quasi-invariant finite measures on the boundary of the Stone-{\v C}ech compactification of $G(\zz)$. In order to do so, we find it convenient to split the measures on this boundary in different classes and treat each type individually. Moreover, it follows from the previous section that we only need to consider measures supported on $\phi^{-1} (M_{3,2})$.
   
  \begin{defn}
 Let $V$ be a $1$-dimensional linear subspace of $\mathbb{R}^3$; we say that $V$ is a \textit{rational line} if it contains a vector whose coordinates are all rational numbers. Let now $W$ be a $2$-dimensional subspace of $\mathbb{R}^3$; we say that $W$ is a \textit{rational plane} if $W^\perp$ is a rational line. 
 \end{defn}

  Let $V, W$ be linear subspaces of $\mathbb{R}^3$ and define $C_{V, W}:=\{(m_1,m_2) \in M_{3,2} \; | \; \ran(m_1)   \subset \mathbb{P}(V), \quad  \mathbb{P}(W) \subset \ker(m_1)\} \subset M_{3,2}$, $D_{V,W} = \{ (m_1, m_2) \in M_{3,2} \; | \; \mathbb{P}(V) \subset \ran (m_1), \quad \ker(m_1) \subset \mathbb{P}(W)\} \subset M_{3,2}$, $E_{V,W} = \{ (m_1,m_2) \in M_{3,2} \; | \;  \mathbb{P}(V) \subset \ran (m_1), \quad \ker (m_1) = \mathbb{P}(W)\}$ and $F_{V,W} = \{(m_1,m_2) \in M_{3,2} \; | \; \ran (m_1)= \mathbb{P}(V), \quad \ker (m_1) \subset \mathbb{P}(W)\}$. We let $\mathcal{C}$ denote the class of sets of the form $C_{V,W}$ such that $0 < \dim (V), \dim (W) < 3$ (i.e. $\dim(V)=2$, $\dim (W) =1$). Similarly, $\mathcal{D}$ is the class of sets of the form $D_{W,V}$, where $\dim (W)=1$, $\dim (V)=2$, $\mathcal{E}$ is the class of sets of the $E_{V,W}$ with $\dim (V) =\dim (W)=1$, $\mathcal{F}$ is the class of sets of the form $F_{V,W}$ with $\dim (V)=\dim (W)=2$. 
  \begin{defn}
Let $V$ and $W$ be linear subspaces of $\mathbb{R}^3$ with $\dim (V)=2$, $\dim (W)=1$. The set $C_{V,W} \in \mathcal{C}$ is called a \textit{(plane, line)-set} and the set $D_{W,V}$ is a \textit{(line, plane)-set}. The (plane, line)-set $C_{V,W}$ is \textit{rational} if $V$ is a rational plane and $W$ is a rational line; it is \textit{irrational} if at least one between $V$ and $W$ is irrational. The sets in $\mathcal{E}$ are called {\it (line, line)-sets} and the sets in $\mathcal{F}$ are called {\it (plane, plane)-sets}. A set which belongs to $\mathcal{D} \cup \mathcal{E} \cup \mathcal{F} \cup \mathcal{C}$ is called a \textit{linear set}. The union of the inverse image under $\phi$ of all the rational (plane, line)- sets is the {\it rational part} (of $\partial_\beta G(\zz)$) and the union of the inverse image under $\phi$ of the orbit of an irrational (plane, line)-set is an {\it irrational part} (of $\partial_\beta G(\zz)$). A Radon measure supported on the rational part is called a {\it rational measure}; a Radon measure supported on an irrational part is called an {\it irrational measure}. A Radon measure on $\phi^{-1} (M_{3,2})$ which annihilates the preimage of every linear set is {\it top-dimensional}.
  \end{defn}


 
 Let now $V, V', W, W'\subset \mathbb{R}^3$ be linear subspaces with $\dim (V)=\dim (V')=1$, $\dim (W)=\dim (W')=2$ and denote by $\braket{V,V'}$ the linear span of $V$ and $V'$. We have the the following possibilities:
 
 \begin{equation*}
 D_{V, W} \cap D_{V', W'} = \begin{cases} D_{V,W} = D_{V',W'} & \mbox{ for } V=V', \; W=W',\\
 E_{V, W \cap W'} =E_{V', W \cap W'}& \mbox{ for } V=V', \; W \neq W',\\
 F_{\braket{V,V'}, W}=F_{\braket{V,V'}, W'} & \mbox{ for } V\neq V', \; W=W',\\
 C_{\braket{V,V'}, W \cap W'} & \mbox{ for } V\neq V', \; W\neq W',
 \end{cases}
 \end{equation*}
 
 If now $V, V', W, W'$ are all $1$-dimensional, we have
 
 \begin{equation*}
 E_{V, W} \cap E_{V',W'} = \begin{cases} E_{V,W}=E_{V',W'} & \mbox{ for } V=V',\; W=W',\\
 C_{\braket{V,V'}, W} = C_{\braket{V,V'}, W'} & \mbox{ for } V \neq V', \; W=W',\\
 \emptyset & \mbox{ for } W\neq W'
 \end{cases}
 \end{equation*}
  
 Similarly, if $V, V', W, W'$ are all $2$-dimensional, we have
 \begin{equation*}
 F_{V, W} \cap F_{V', W'} = \begin{cases} F_{V, W} =F_{V', W'} & \mbox{ for } V=V', \; W=W',\\
 C_{V, W \cap W'} = C_{V', W \cap W'} & \mbox{ for } V=V', \; W \neq W',\\
 \emptyset & \mbox{ for } V \neq V'
 \end{cases}
 \end{equation*}
  
  If now $V$ and $V'$ have dimension $2$, while $W$ and $W'$ have dimension $1$, we have
   \begin{equation*}
 C_{V, W} \cap C_{V', W'} = \begin{cases} C_{V, W} =C_{V', W'} & \mbox{ for } V=V', \; W=W',\\
 \emptyset & \mbox{ otherwise }.
 \end{cases}
 \end{equation*}
  
Consider now the set of linear subspaces of $\mathbb{R}^3$ endowed with the natural action of $G(\mathbb{Z})$.

 \begin{lem}
  \label{intersections}
 \label{globalstabsets}
 Let $V,V', W, W' \subset \mathbb{R}^3$ be linear subspaces with $\dim (V)=\dim (V')=2$, $\dim (W)=\dim (W')=1$, we have
 \begin{equation*}
 (G(\mathbb{Z}) \times G(\mathbb{Z}))_{D_{W,V}} = G(\mathbb{Z})_W \times G(\mathbb{Z})_V,
 \end{equation*}
 \begin{equation*}
  (G(\mathbb{Z}) \times G(\mathbb{Z}))_{E_{W,W'}} = G(\mathbb{Z})_W \times G(\mathbb{Z})_{W'},
  \end{equation*}
  \begin{equation*}
   (G(\mathbb{Z}) \times G(\mathbb{Z}))_{F_{V,V'}} = G(\mathbb{Z})_V \times G(\mathbb{Z})_{V'},
   \end{equation*}
   \begin{equation*}
    (G(\mathbb{Z}) \times G(\mathbb{Z}))_{C_{V,W}} = G(\mathbb{Z})_V \times G(\mathbb{Z})_W,
    \end{equation*}
    In particular, the global stabilizers of these types of sets are groups.
   \end{lem}
 \proof Let $(\gamma_1, \gamma_2) \in G(\mathbb{Z}) \times G(\mathbb{Z})$ and suppose $\gamma_1 \notin G(\mathbb{Z})_W$. Then $\gamma_1$ sets up a bijection between the planes containing $W$ and the planes containing $\gamma_1 (W) \neq W$. Since these two sets are different, it follows that there is a plane $\pi$ containing $W$ whose image under $\gamma_1$ does not contain $W$. Hence it is enough to choose a matrix of rank $2$ whose image is $\pi$ and kernel is contained in $V$ in order to see that $(\gamma_1, \gamma_2)$ cannot belong to $(G(\mathbb{Z}) \times G(\mathbb{Z}))_{D_{W,V}}$. In the same way, if $\gamma_2 (V) \neq V$, there is a $1$-dimensional subspace $r$ of $V$ whose image under $\gamma_2$ does not belong to $V$. Hence, for every $m$ with $\ker(m) = r$ we have $\ker (m \gamma_2^{-1}) = \gamma_2 (\ker (m)) \notin V$ and it follows that $(G(\mathbb{Z}) \times G(\mathbb{Z}))_{D_{W,V}} = G(\mathbb{Z})_W \times G(\mathbb{Z})_V$. The other cases are proved in the same way. $\Box$
 
 
\subsection{Top-dimensional measures}
\label{td}
Inspired by \cite{ams} and \cite{GoMa}, in this section we use proximality arguments in order to obtain information on the decay of matrix coefficients for Koopman representations associated to top-dimensional measures.\\
Let $k \in \mathbb{N}$. We denote by $M_{k,\| \cdot \|=1}(\mathbb{R})$ the compact set of $k \times k$ norm one matrices.

 \begin{lem}
 \label{lemams}
 Let $k \in \mathbb{N}$, $\{x_n\}$, $\{y_n\}$ and $\{z_n\}$ be sequences in $M_k (\mathbb{R})$ such that there are norm one matrices $x$, $y$ and $z$ satisfying
 \begin{equation*}
 x_n /\| x_n\| \rightarrow x, \quad y_n /\| y_n\| \rightarrow y, \quad z_n /\|z_n\| \rightarrow z, \quad xyz \neq 0.
 \end{equation*}
  Then
  \begin{itemize}
  \item[(i)]: the sequence $x_n y_n z_n /\| x_n y_n z_n\|$ converges to a point $l \in M_{k, \|\cdot \|=1}$,
  \item[(ii)]: there is $\lambda \in \mathbb{R}_+ \backslash \{0\}$ such that $l = \lambda xyz$.
  \end{itemize}
  \end{lem}
  \proof This follows from the proof of \cite{ams} Lemma 4.3 (b). We include the details for completeness. By compactness of $M_{k, \| \cdot \|=1}(\mathbb{R})$ it is enough to show that there is $0 \neq \lambda \in \mathbb{R}_+$ such that every convergent subsequence of $x_n y_n z_n /\| x_n y_n z_n\|$ converges to $\lambda xyz$. Up to reindexing we can then suppose $x_n y_n z_n /\| x_n y_n z_n\| \rightarrow c \in M_{k, \|\cdot \|=1}$. Let $v \in \mathbb{R}^k$ be such that $xy z(v) \neq 0$. We have
  \begin{equation}
  \label{eqc}
  \frac{x_n y_n z_n}{\| x_n \| \| y_n\| \|z_n\|} (v) \rightarrow xyz (v), \qquad \frac{x_n y_n z_n}{\| x_n y_n z_n\|} (v) \rightarrow c(v) \neq 0 \quad \mbox{(since $\| c(v) \| \geq \| xy z(v)\|$)}.
  \end{equation}
  Denoting $\lambda_n = \| x_n y_n z_n \| / (\|x_n \| \| y_n\| \|z_n\|)$ we can rewrite (\ref{eqc}) as
  \begin{equation*}
  \frac{x_n y_n z_n}{\|x_n y_n z_n\|} (v) \rightarrow c(v) \neq 0, \qquad  \lambda_n \frac{x_n y_n z_n}{\|x_n y_n z_n\|}(v) \rightarrow xyz(v) \neq 0.
  \end{equation*}
  We want to check that under this condition $\lambda_n$ converges to a non-zero real number $\lambda$. Indeed, suppose this is not the case, then $\lambda_n c(v)$ would not converge, but $|\lambda_n   | \leq 1$ for every $n$ and so
  \begin{equation*}
  \|\lambda_n c(v) - xyz(v)\| \leq \|\lambda_n c(v) - \lambda_n \frac{x_n y_n z_n}{\| x_n y_n z_n\|} (v) \| + \| \lambda_n \frac{x_n y_n z_n}{\|x_n y_n z_n\|} (v) - xyz(v)\| \rightarrow 0.
  \end{equation*}
  Then there is $0 \neq \lambda \in \mathbb{R}_+$ such that $\lambda_n \rightarrow \lambda$ and so
  \begin{equation*}  
  \begin{split}
  \left\|   \frac{x_n y_n z_n}{\| x_n y_n z_n\|}  - \lambda^{-1} xyz \right\|&  \leq \left\|   \frac{x_n y_n z_n}{\| x_n y_n z_n\|} - \lambda^{-1}  \frac{x_n y_n z_n}{\| x_n \| \| y_n\| \|z_n\|}\right\| + \left\| \lambda^{-1}  \frac{x_n y_n z_n}{\| x_n \| \| y_n\| \|z_n\|} - \lambda^{-1} xyz \right\|    \\
  & \leq   \left\|   \frac{x_n y_n z_n}{\| x_n y_n z_n\|} - \lambda^{-1} \lambda_n \frac{x_n y_n z_n}{\| x_n   y_n z_n\|}\right\| + \left\| \lambda^{-1}  \frac{x_n y_n z_n}{\| x_n \| \| y_n\| \|z_n\|} - \lambda^{-1} xyz \right\|;  
  \end{split}
  \end{equation*}
this inequality shows that $xyz=\lambda c$. $\Box$

In the next Lemma we will use the following trivial consequence of Lemma \ref{lemams}: if $x_n$, $y_n$ and $z_n$ are sequences in $M_3 (\mathbb{R})$ converging respectively to $x$, $y$ and $z$ and $xyz \neq 0$, then $\mathbb{P} (x_n y_n z_n) \rightarrow \mathbb{P}(xyz)$.     
     \begin{lem}
     \label{lemmapl}
     Let $g_n = (\gamma_n , \eta_n)$ be a sequence in $G(\mathbb{Z}) \times G(\mathbb{Z})$ such that both $\gamma_n$ and $\eta_n$ diverge. Then there are a subsequence $g_{n_k} = (\gamma_{n_k}, \eta_{n_k})$ and (line-plane)-sets $D^0_1,...,D^0_m$, $D_1,...,D_l$, $m,l \in \ \mathbb{N}$, such that for every $y \in M_{3,2} \backslash \bigcup_{i=1}^m D^0_i$  we have
     \begin{itemize}
     \item[(i)] the sequence $(\gamma_{n_k}, \eta_{n_k})\cdot y$ converges;
     \item[(ii)] the limit $\lim_k (\gamma_{n_k}, \eta_{n_k}) \cdot y$ either has rank $1$ or belongs to $\bigcup_{i=1}^l D_i$.
     \end{itemize}
    \end{lem}
     \proof  It is enough to find sets $D_i$ in $\mathcal{C} \cup \mathcal{D} \cup \mathcal{E} \cup \mathcal{F}$ rather than just in $\mathcal{D}$ with the above properties in order to prove the statement. Up to taking a subsequence we can suppose that $\gamma_n /\| \gamma_n\|$ and $\eta_n / \| \eta_n\|$ converge to $\gamma=\gamma^{(0)}$ and $\eta=\eta^{(0)}$ respectively. Given a linear space $L \subset \mathbb{R}^3$ denote by $P_L$ the orthogonal projection with range $L$. Let then $\gamma_n^{(1)} := \gamma_n P_{\ker (\gamma)} $, $\eta_n^{(1)} := P_{\ran (\eta)^\perp} \eta_n $. Again, up to taking a subsequence, we can suppose that both $\gamma_n^{(1)} /\| \gamma_n^{(1)}\|$ and $\eta_n^{(1)} /\| \eta_n^{(1)}\|$ converge respectively to $\gamma^{(1)}$ and $\eta^{(1)}$. Since $\ker(\gamma)^\perp  \subset \ker (\gamma^{(1)})$ we have $\ker(\gamma^{(1)})^\perp \subset \ker(\gamma)$; similarly, since $\ran(\eta^{(1)}) \subset \ran (\eta)^\perp$, we have $\ran (\eta )  \subset \ran (\eta^{(1)})^\perp $. If $\ker (\gamma^{(1)})^\perp \neq \ker (\gamma)$, we have $\ker (\gamma) \cap \ker (\gamma^{(1)}) \neq \{0\}$ and we consider $\gamma_n^{(2)} := \gamma_n P_{\ker (\gamma) \cap \ker (\gamma^{(1)})}$; up to taking a subsequence, we have $\gamma_n^{(2)} /\| \gamma_n^{(2)}\| \rightarrow \gamma^{(2)}$. If $\ran (\eta) \neq \ran (\eta^{(1)})^\perp$, then $\ran (\eta^{(1)})^\perp \cap \ran (\eta)^\perp \neq \{0\}$ and we consider $\eta_n^{(2)} :=P_{\ran(\eta)^\perp \cap \ran(\eta^{(1)})^\perp} \eta_n $; up to taking a subsequence we have $\eta_n^{(2)} / \| \eta_n^{(2)}\| \rightarrow \eta^{(2)}$. Again, up to taking a subsequence, we can suppose that the following sequences converge in $\mathbb{R} \cup \{\infty\}$: 
     \begin{equation*}
     \frac{\| \gamma_n P_{\ker (\g)} \|}{\| \gamma_n\|} \frac{\|\eta_n\|}{\| P_{\ran (\e)^\perp} \eta_n\|}, \quad \frac{\|\gamma_n \|}{\| \gamma_n P_{\ker (\g)}\|} \frac{\| P_{\ran (\e)^\perp \cap \ran(\euno)^\perp} \eta_n \|}{\| \eta_n\|},
     \end{equation*}
     \begin{equation*}
     \frac{\| \eta_n \|}{ \| P_{\ran (\e)^\perp} \eta_n \|}\frac{ \| \gamma_n P_{\ker (\g) \cap \ker (\guno)} \|}{ \| \gamma_n \|}, \quad  \frac{\| \eta_n \|}{ \| P_{\ran (\e)^\perp} \eta_n \|}\frac{ \| \gamma_n P_{\ker (\g) \cap \ker (\guno)} \|}{ \| \gamma_n P_{\ker (\g)}\|},
     \end{equation*}
     \begin{equation*}
     \frac{\|\gamma_n\|}{\| \gamma_n P_{\ker (\g)}\|} \frac{ \| P_{\ran (\e)^\perp \cap \ran (\euno)^\perp} \eta_n \|}{\| P_{\ran (\e)^\perp} \eta_n\|}.
     \end{equation*}
     Note that the following relations hold:
     \begin{equation*}
     \frac{\|\gamma_n P_{\ker (\gamma^{(0)})}\|}{\|\gamma_n\|} \rightarrow 0, \quad \frac{\|\gamma_n P_{\ker (\gamma^{(0)} )\cap \ker(\gamma^{(1)})}\|}{\|\gamma_n P_{\ker (\gamma^{(0)})}\|} \rightarrow 0,
     \end{equation*}
     \begin{equation*}
     \frac{\|P_{\ran (\eta^{(0)})^\perp \eta_n \|}}{\| \eta_n\|} \rightarrow 0, \quad \frac{\|P_{\ran (\eta^{(0)})^\perp \cap \ran (\eta^{(1)})^\perp} \eta_n \|}{\|P_{ \ran(\eta^{(0)})^\perp} \eta_n\|} \rightarrow 0
     \end{equation*}
     and that we have the following possibilities for the ranks of these matrices (we assume $\gdue =0$ if $\ker (\gamma^{(1)})^\perp= \ker (\gamma)$ and $\edue =0$ if $\ran (\eta) =\ran (\eta^{(1)})^\perp$):
     \begin{equation*}
     (\rank (\g), \rank(\guno), \rank(\gdue)) \in \{ (2,1,0), (1,2,0), (1,1,1)\},
     \end{equation*}
      \begin{equation*}
     (\rank (\e), \rank(\euno), \rank(\edue)) \in \{ (2,1,0), (1,2,0), (1,1,1)\}.
     \end{equation*}
     Suppose now that $x$ is a $3 \times 3$-matrix of rank $2$ and norm $1$, which is a representative for $m_1$, where $y=(m_1 , m_2) \in M_{3,2}$. We distinguish several different cases.
     
$(i)$:     Suppose $\gamma_0 x \eta_0 \neq 0$. Then $\gamma^{(0)} x \eta^{(0)} = \lim_n (\gamma_n x \eta_n)/( \|\gamma_n\|  \|\eta_n\|)$ and so the limit $\mathbb{P} (\gamma_n x \eta_n)$ either has rank $1$ (if $(\rank (\gamma^{(0)}), \rank (\eta^{(0)})) \neq (2,2)$) or belongs to $C_{\ran(\gamma^{(0)}), \ker(\eta^{(0)})}$.
     
 $(ii)$:    Suppose now that
     \begin{equation*}
    \begin{array}{cc} \g x \e =0, & \g x \euno =0\\
    \guno x \e =0, & \guno x \euno=0.
    \end{array}
    \end{equation*}
     First of all observe that in this istance we should have $\rank(\e)=\rank(\euno)=\rank(\g)=\rank(\guno)=1$. Indeed if one between $\e$ and $\euno$ has rank $2$, the condition implies $\ran (x) \subset \ker(\g) \cap \ker(\guno)$, which is at most one-dimensional; on the other hand, if one between $\g$ and $\guno$ has rank $2$, then $\ker (\g) \cap \ker (\guno)=\{0\}$ and so the kernel of $x$ should contain the linear space generated by $\ran (\e)$ and $\ran(\euno)$, contradicting the fact that $x$ has rank $2$. Now, if $x(\braket{\ran (\e), \ran(\euno)}) = \ker(\g) \cap \ker(\guno)$, then $\ker (x) \subset \braket{\ran(\e), \ran(\euno)}$ and so $m$ belongs to $D_{\ker (\g) \cap \ker (\guno), \braket{\ran(\e), \ran(\euno)}}$.
     
     $(iii)$: Consider the case in which
      \begin{equation*}
    \begin{array}{cc} \g x \e =0, & \g x \euno \neq 0\\
    \guno x \e \neq 0. & \empty
    \end{array}
    \end{equation*}
     We have
     \begin{equation*}
     \guno x \e = \guno x P_{\ran (\e)} \e = \lim \frac{\gamma_n P_{\ker (\g)} x P_{\ran (\e)} \eta_n}{\| \gamma_n P_{\ker (\g) }\| \| \eta_n\|},
     \end{equation*}
     \begin{equation*}
     \g x \euno = \lim \frac{\gamma_n x P_{\ran (\e)^\perp} \eta_n}{\|\gamma_n \| \|P_{\ran (\e)^\perp }\eta_n\|}.
     \end{equation*}
 Since $x P_{\ran (\e)} = P_{\ker (\g)} x P_{\ran (\e)}$, we can write
 \begin{equation*}
 \begin{split}
 \frac{\gamma_n x \eta_n}{\|\gamma_n \| \|P_{\ran (\e)^\perp} \eta_n\|} &=  \frac{\gamma_n P_{\ker (\g)}x P_{\ran (\e)}\eta_n}{\|\gamma_n \| \|P_{\ran (\e)^\perp} \eta_n\|} +  \frac{\gamma_n x P_{\ran (\e)^\perp} \eta_n}{\|\gamma_n \| \|P_{\ran (\e)^\perp} \eta_n\|}\\
 &=\left( \frac{\|\gamma_n P_{\ker (\g)} \|}{\| P_{ \ran (\e)^\perp} \eta_n \|} \cdot \frac{ \|\eta_n\|}{\|\gamma_n\|} \right) \frac{\gamma_n P_{\ker (\g)} x P_{\ran (\e)} \eta_n}{\| \gamma_n P_{\ker (\g) }\|\| \eta_n\|} \\ &+  \frac{\gamma_n x P_{\ran (\e)^\perp} \eta_n}{\|\gamma_n \| \|P_{\ran (\e)^\perp} \eta_n\|}.
 \end{split}
 \end{equation*}
     Hence, letting $\lambda_n = \left(\frac{\|\gamma_n P_{\ker (\g)} \|}{\| P_{ \ran (\e)^\perp} \eta_n \|} \cdot \frac{ \|\eta_n\|}{\|\gamma_n\|} \right)$, we have the following possibilities
     \begin{equation*}
     \lim_n \mathbb{P}(\gamma_n x \eta_n)= \begin{cases} \mathbb{P}( \guno x \e) & \mbox{ for } \lambda_n \rightarrow \infty\\
     \mathbb{P}(\g x \euno) & \mbox{ for } \lambda_n \rightarrow 0\\
     \mathbb{P} (\lambda \guno x \e + \g x \euno) & \mbox{ for } \lambda_n \rightarrow \lambda \in \mathbb{R} \backslash \{0\}
     \end{cases}.
     \end{equation*}
     Hence, in the first two cases, if the limit has rank $2$, then it belongs to either $C_{\ran (\guno), \ker \e}$ or $C_{\ran (\g), \ker (\euno)}$. Consider the third situation. If both $\g x \euno$ and $\guno x \e$ have rank $1$, then their ranges should be linearly independent (if the limit does not have rank $1$) and so the kernel of $\lambda \guno x \e + \g x \euno$ is the intersection of their kernels; hence, if $i$ is such that $\rank (\gamma^{(i)})=1$ and $j$ is such that $\rank (\eta^{(j)})=1$, then $  \lim_n \mathbb{P}(\gamma_n x \eta_n)$ belongs to $D_{\ran (\gamma^{(i)}), \ker (\eta^{(j)})}$. This is clearly the case if $\rank (\gamma^{(i)})= \rank (\eta^{(i)})=1$ for every $i \in \{0,1\}$. Suppose that $\rank(\g)=\rank (\euno)=1$ and $\rank (\guno)=\rank (\e)=2$. We have $x\ran (\e) \subset \ker (\g)$, now the case $\ran (x\e)= \ker (\g)$ is impossible, since this would entail $m \ran (\euno) \subset \ker (\g)$, which is in contrast with the assumption $\g x \euno \neq 0$. Hence the image under $x$ of $\ran (\e)$ is a $1$-dimensional subspace of $\ker (\g)$ and so $\guno x \e$ has rank $1$. Suppose now that $\rank(\g)=\rank (\euno)=2$ and $\rank (\guno)=\rank (\e)=1$. Then we can reduce to the above situation by considering ${}^t x$, obtaining $\mathbb{P}(\lambda {}^t \e {}^t x {{}^t\guno} + {}^t \euno {}^t x {}^t \g) \in D_{\ran ({}^t \eta^{(i)}), \ker ({}^t \gamma^{(j)})}$ for some $i,j$. Hence $\mathbb{P}(\lambda  \guno x {\e} + \g x \euno)$ belongs to $D_{ \ker ({}^t \gamma^{(j)})^\perp,\ran ({}^t \eta^{(i)})^\perp}= D_{\ran (\gamma^{(j)}), \ker (\eta^{(i)})}$.
     
$(iv)$: Consider now the case
\begin{equation*}
  \begin{array}{cc} \g x \e =0, & \g x \euno = 0\\
    \guno x \e \neq 0. & \empty
    \end{array}
    \end{equation*}    
In order to deal with this case we consider several subcases. Suppose
\begin{equation*}
\g x \edue =0.
\end{equation*}       
Hence $\ran x \subset \ker (\gamma^{(0)})$ and so
\begin{equation*}
\frac{\gamma_n x \eta_n}{\| \gamma_n P_{\ker (\gamma^{(0)})} \| \| \eta_n\|} = \frac{\gamma_n P_{\ker (\g)} x \eta_n}{\| \gamma_n P_{\ker (\g)} \| \| \eta_n\|} \rightarrow \guno x \e.
\end{equation*}     
Suppose now that
\begin{equation*}
\g x \edue \neq 0, \qquad \guno x \euno \neq 0.
\end{equation*}
Then we can write
\begin{equation*}
\begin{split}
\frac{\gamma_n x \eta_n}{\| \gamma_n P_{\ker (\g)} \| \|\eta_n\|} &= \frac{\gamma_n P_{\ker(\g)} x P_{\ran (\e)} \eta_n + \gamma_n x P_{\ran (\e)^\perp} \eta_n}{\| \gamma_n P_{\ker (\g)} \|  \|\eta_n\|} \\
& = \frac{\gamma_n P_{\ker(\g)} x P_{\ran (\e)} \eta_n}{\| \gamma_n P_{\ker (\g)} \|  \|\eta_n\|} +\frac {\gamma_n x P_{\ran (\euno) \cap \ran (\e)^\perp} \eta_n}{\| \gamma_n P_{\ker (\g)} \|  \|\eta_n\|} \\ & + \frac{\gamma_n x P_{\ran(\euno)^\perp \cap \ran(\e)^\perp} \eta_n}{\| \gamma_n P_{\ker (\g)} \|  \|\eta_n\|}\\
& = \frac{\gamma_n P_{\ker(\g)} x P_{\ran (\e)} \eta_n}{\| \gamma_n P_{\ker (\g)} \|  \|\eta_n\|} \\
& + \left( \frac{\|P_{\ran(\e)^\perp} \eta_n\|}{\| \eta_n\|} \right) \frac{\gamma_n P_{\ker (\g)} x P_{\ran (\euno)} \eta_n}{\| \gamma_n P_{\ker (\g)} \|  \| P_{\ran (\e)^\perp} \eta_n\|}\\
& +  \left( \frac{\|\gamma_n\| \| P_{\ran (\e)^\perp \cap \ran(\euno)^\perp} \eta_n \|}{ \| \gamma_n P_{\ker (\g)} \| \|\eta_n\|} \right) \frac{ \gamma_n x P_{\ran (\e)^\perp \cap \ran (\euno)^\perp} \eta_n}{\|\gamma_n\| \|P_{\ran (\e)^\perp \cap \ran (\euno)^\perp} \eta_n\|}.
\end{split}
\end{equation*}
Since we have
     \begin{equation*}
 \frac{\gamma_n P_{\ker(\g)} x P_{\ran (\e)} \eta_n}{\| \gamma_n P_{\ker (\g)}\| \|\eta_n\|}  \rightarrow \guno x \e,
 \end{equation*}
 \begin{equation*}
  \frac{\gamma_n P_{\ker (\g)} x P_{\ran (\euno)} \eta_n}{\| \gamma_n P_{\ker (\g)} \| \| P_{\ran (\e)^\perp} \eta_n\|}\rightarrow \guno x \euno,
  \end{equation*}
  \begin{equation*}
     \frac{ \gamma_n x P_{\ran (\e)^\perp \cap \ran (\euno)^\perp} \eta_n}{\|\gamma_n\| \|P_{\ran (\e)^\perp \cap \ran (\euno)^\perp}\eta_n \|} \rightarrow \g x \edue,
      \end{equation*}
      \begin{equation*}
       \frac{\|P_{\ran(\e)^\perp} \eta_n\|}{\| \eta_n\|}  \rightarrow 0,
       \end{equation*}
      we have the following possibilities, depending on the value of the limit $\lambda$ of the sequence $\lambda_n =\frac{\|\gamma_n\| \| P_{\ran (\e)^\perp \cap \ran (\euno)^\perp} \eta_n\|}{\| \gamma_n P_{ \ker (\g)}\| \| \eta_n\|}$:
      \begin{equation*}
      \lim_n \mathbb{P} (\gamma_n x \eta_n) = \begin{cases}	\mathbb{P}(\guno x \e) & \mbox{ for } \lambda = 0\\
      \mathbb{P} (\g x \edue) & \mbox{ for } \lambda =\infty\\
      \mathbb{P}(\guno x \e + \lambda \g x \edue)& \mbox{ for } \lambda \neq 0, \infty.
      \end{cases}
     \end{equation*}
In the first two cases the result belongs to a set of the form $C_{\ran (\gamma^{(i)}), \ker (\eta^{(j)})}$, while in the third case both $\guno x \e$ and $\g x \edue$ have rank $1$; hence $\lim_n \mathbb{P}(\gamma_n x \eta_n)$ belongs to $D_{\ran (\gamma^{(i)}), \ker (\eta^{(j)})}$ for some $i,j$ with $\rank (\gamma^{(i)})=\rank (\eta^{(j)})=1$.     

Let now
\begin{equation*}
\begin{array}{ccc} \g x \e =0, & \g x \euno = 0, & \g x \edue \neq 0\\
    \guno x \e \neq 0, & \guno x \euno =0, & \\
    \gdue x \euno \neq 0.
    & \empty &
    \end{array}
\end{equation*}     
We write
\begin{equation*}
\begin{split}
\frac{\gamma_n x \eta_n}{\|\gamma_n P_{\ker (\g)} \|\| \eta_n\|}    & =
\frac{\gamma_n x (P_{\ran (\e)} + P_{\ran(\e)^\perp}) \eta_n}{\|\gamma_n P_{\ker (\g)} \|\| \eta_n\|}\\ 
& = \frac{\gamma_n P_{\ker (\g)}x P_{\ran (\e)}\eta_n}{\|\gamma_n P_{\ker (\g)} \|   \| \eta_n\|} + \frac{\gamma_n x P_{\ran (\e)^\perp} \eta_n}{\|\gamma_n P_{\ker (\g)} \|   \| \eta_n\|}\\
& = \frac{\gamma_n P_{\ker (\g)}x P_{\ran (\e)}\eta_n}{\|\gamma_n P_{\ker (\g)} \|   \| \eta_n\|}\\
& + \frac{\gamma_n x P_{\ran (\e)^\perp \cap \ran (\euno)}\eta_n}{\|\gamma_n P_{\ker (\g)} \|   \| \eta_n\|} + \frac{\gamma_n x P_{\ran (\e)^\perp \cap \ran (\euno)^\perp}\eta_n}{\|\gamma_n P_{\ker (\g)} \|   \| \eta_n\|} \\
& =  \frac{\gamma_n P_{\ker (\g)}x P_{\ran (\e)}\eta_n}{\|\gamma_n P_{\ker (\g)} \|   \| \eta_n\|}\\
& + \frac{\gamma_n P_{\ker (\g) \cap \ker (\guno)}x P_{\ran (\e)^\perp \cap \ran (\euno)}\eta_n}{\|\gamma_n P_{\ker (\g)} \|   \| \eta_n\|} \\
&  + \frac{\gamma_n x P_{\ran (\e)^\perp \cap \ran (\euno)^\perp}\eta_n}{\|\gamma_n P_{\ker (\g)} \|   \| \eta_n\|} \\
& =   \frac{\gamma_n P_{\ker (\g)}x P_{\ran (\e)}\eta_n}{\|\gamma_n P_{\ker (\g)} \|  \| \eta_n\|}\\
&  + \left( \frac{\| P_{\ran (\e)^\perp} \eta_n\| \| \gamma_n P_{\ker (\g) \cap \ker (\guno)} \|}{\| \eta_n \| \| \gamma_n P_{\ker (\g)}\|} \right) \frac{\gamma_n P_{\ker (\g) \cap \ker (\guno)}x P_{\ran (\euno)} P_{\ran(\e)^\perp}\eta_n}{\|\gamma_n P_{\ker (\g) \cap \ker (\guno)} \|   \|P_{\ran(\e)^\perp} \eta_n\|} \\
& + \left( \frac{ \| \gamma_n\| \| P_{\ran (\e)^\perp \cap \ran (\euno)^\perp} \eta_n \|}{\| \gamma_n P_{\ker (\g)}\|\| \eta_n\|} \right) \frac{\gamma_n x P_{\ran (\e)^\perp \cap \ran (\euno)^\perp}\eta_n}{\|\gamma_n\| \| P_{\ran(\e)^\perp \cap \ran (\euno)^\perp}\eta_n\|} .
 \end{split}
 \end{equation*}    
     Since
     \begin{equation*}
     \lim_n \left( \frac{\| P_{\ran (\e)^\perp} \eta_n\| \| \gamma_n P_{\ker (\g) \cap \ker (\guno)} \|}{\| \eta_n \| \| \gamma_n P_{\ker (\g)}\|} \right) =0,
     \end{equation*}
 denoting by $\lambda$ the limit of $\frac{\|\gamma_n \| \|P_{\ran (\e)^\perp \cap \ran (\euno)^\perp} \eta_n \|}{\| \gamma_n P_{\ker (\g)} \| \| \eta_n\|}$, we obtain   
 \begin{equation*}
\lim_n \mathbb{P} (\gamma_n x \eta_n) = \begin{cases} \mathbb{P}(\guno x \e) & \mbox{ for } \lambda =0\\
\mathbb{P}( \g x \edue) & \mbox{ for } \lambda = \infty\\
\mathbb{P} (\guno x \e + \lambda \g x \edue) & \mbox{ for } \lambda \neq 0, \infty 
\end{cases}.
\end{equation*}
    Hence we are in the same situation of the previous subcase.
    
  Let now  
    \begin{equation*}
    \begin{array}{ccc} \g x \e =0, & \g x \euno = 0, & \g x \edue \neq 0\\
    \guno x \e \neq 0, & \guno x \euno =0, & \\
    \gdue x \euno = 0.
    & \empty &
    \end{array}
\end{equation*}     
    In this case $\ran (\euno) \subset \ker (x)$ and so we can write
    \begin{equation*}
    \begin{split}
    \frac{\gamma_n x \eta_n}{\| \gamma_n P_{\ker (\g)}\| \|\eta_n\|}& =
    \frac{\gamma_n x P_{\ran (\e)} \eta_n + \gamma_n x P_{\ran (\e)^\perp \cap \ran(\euno)^\perp} \eta_n + \gamma_n x P_{\ran (\e)^\perp \cap \ran (\euno)} \eta_n}{\| \gamma_n P_{\ker (\g)}\|    \|\eta_n\|}\\
    &  =
    \frac{\gamma_n x P_{\ran (\e)} \eta_n + \gamma_n x P_{\ran (\e)^\perp \cap \ran(\euno)^\perp} \eta_n }{\| \gamma_n P_{\ker (\g)}\|    \|\eta_n\|}\\
    & = \frac{\gamma_n P_{\ker (\g)} x P_{\ran (\e)} \eta_n}{\| \gamma_n P_{ \ker (\g)} \|  \| \eta_n\|}\\
    & +  \left (\frac{\|\gamma_n \| \| P_{\ran (\e)^\perp \cap \ran (\euno)^\perp} \eta_n \|}{ \| \gamma_n P_{\ker (\g)}\| \| \eta_n\|} \right) \frac{\gamma_n x P_{\ran (\e)^\perp \cap \ran (\euno)^\perp} \eta_n}{\| \gamma_n \| \| P_{\ran (\e)^\perp \cap \ran (\euno)^\perp} \eta_n \|}.
     \end{split}
     \end{equation*}
Letting $\lambda = \lim_n \frac{\| \gamma_n \| \| P_{\ran (\e)^\perp \cap \ran (\euno)^\perp} \eta_n\|}{\| \gamma_n P_{\ker (\g)}\| \| \eta_n\|}$, we obtain    
      \begin{equation*}
\lim_n \mathbb{P} (\gamma_n x \eta_n) = \begin{cases} \mathbb{P}(\guno x \e) & \mbox{ for } \lambda =0\\
\mathbb{P}( \g x \edue) & \mbox{ for } \lambda = \infty\\
\mathbb{P} ( \guno x \e + \lambda \g x \edue) & \mbox{ for } \lambda \neq 0, \infty 
\end{cases}.
\end{equation*}
This situation is thus equivalent to the ones considered above.
     
$(v)$: Consider now the case
\begin{equation*}
\begin{array}{cc}	\g x \e =0, & \g x \euno \neq 0,\\
	\guno x \e = 0. & 
	\end{array}
	\end{equation*}
     As above, we will study different subcases. First, we add the condition
     \begin{equation*}
     \gdue x \e \neq 0.
     \end{equation*}
     We have
     \begin{equation*}
     \begin{split}
     \frac{\gamma_n x \eta_n}{\| \gamma_n P_{\ker (\g) \cap \ker (\euno)} \|  \| \eta_n\|} &= \frac{\gamma_n x P_{\ran (\e)} \eta_n + \gamma_n x P_{\ran (\e)^\perp} \eta_n}{\| \gamma_n P_{\ker (\g) \cap \ker (\euno)} \|     \| \eta_n\|}\\
   &  = \frac{\gamma_n P_{ \ker (\g) \cap \ker (\guno)}x P_{\ran (\e)} \eta_n + \gamma_n x P_{\ran (\e)^\perp} \eta_n}{\| \gamma_n P_{\ker (\g) \cap \ker (\euno)} \|     \| \eta_n\|}\\
   & =\frac{\gamma_n P_{ \ker (\g) \cap \ker (\guno)}x P_{\ran (\e)} \eta_n}{\| \gamma_n P_{\ker (\g) \cap \ker (\euno)} \|     \| \eta_n\|}\\
   & + \left( \frac{\| \gamma_n\| \| P_{\ran( \e)^\perp} \eta_n \|}{ \| \gamma_n P_{\ker (\g) \cap \ker (\guno)}\| \| \eta_n \|} \right) \frac{\gamma_n x P_{\ran (\e)^\perp} \eta_n}{\| \gamma_n \| \| P_{\ran (\e)^\perp} \eta_n \|}.
   \end{split}
   \end{equation*}
  Hence, depending on the value of the limit $\lambda$ of $\frac{\|\gamma_n \| \| P_{\ran(\e)^\perp} \eta_n \|}{ \| \gamma_n P_{\ker (\g) \cap \ker (\guno)} \| \| \eta_n\| }$, we obtain   
  
      \begin{equation*}
\lim_n \mathbb{P} (\gamma_n x \eta_n) = \begin{cases} \mathbb{P}(\gdue x \e) & \mbox{ for } \lambda =0\\
\mathbb{P}( \g x \euno) & \mbox{ for } \lambda = \infty\\
\mathbb{P} ( \gdue x \e + \lambda \g x \euno) & \mbox{ for } \lambda \neq 0, \infty 
\end{cases}.
\end{equation*}
Thus, if the limit has rank $2$, it belongs to a set of the form $D_{\ran (\gamma^{(i)}), \ker (\eta^{(j)})}$ with $\rank (\gamma^{(i)})= \rank (\eta^{(j)})=1$.
  
  Suppose now that 
  \begin{equation*}
  \begin{array}{cc}	\g x \e =0, & \g x \euno \neq 0,\\
	\guno x \e = 0, &  \gdue x \e =0.
	\end{array}
	\end{equation*}
In this case $\ran (\e) \subset \ker (x)$ and so
\begin{equation*}
\frac{\gamma_n x \eta_n}{ \| \gamma_n \| \| P_{\ran (\e)^\perp } \eta_n \| } = \frac{\gamma_n x P_{\ran (\e)^\perp} \eta_n}{ \| \gamma_n \| \| P_{\ran (\e)^\perp \eta_n}\|} \rightarrow \g x \euno.
\end{equation*}  
  
  $(vi)$: Let now
  \begin{equation*}
  \begin{array}{cc}
  \g x \e =0, & \g x \euno =0,\\
  \guno x \e =0, & \guno x \euno \neq 0.
  \end{array}
  \end{equation*}
  Suppose first that we also have
  \begin{equation*}
  \g x \edue =0, \quad \gdue x \e \neq 0.
  \end{equation*}
  In this case we have $\ran (x) \subset \ker (\g)$ and we can write
  \begin{equation*}
  \begin{split}
  \frac{\gamma_n x \eta_n}{\| \gamma_n P_{\ker (\g) \cap \ker (\guno)} \|  \| \eta_n\|}&=  \frac{\gamma_n P_{\ker (\g)} x \eta_n}{\| \gamma_n P_{\ker (\g) \cap \ker (\guno)} \|     \| \eta_n\|}\\
  &=  \frac{\gamma_n P_{\ker (\g)} x P_{\ran (\e)}\eta_n}{\| \gamma_n P_{\ker (\g) \cap \ker (\guno)} \|     \| \eta_n\|}\\
  &    +  \frac{\gamma_n P_{\ker (\g)} x P_{\ran (\e)^\perp}\eta_n}{\| \gamma_n P_{\ker (\g) \cap \ker (\guno)} \|     \| \eta_n\|}\\
  &    =  \frac{\gamma_n P_{\ker (\g) \cap \ker (\guno)} x P_{\ran (\e)}\eta_n}{\| \gamma_n P_{\ker (\g) \cap \ker (\guno)} \|     \| \eta_n\|}\\
  &    +  \frac{\gamma_n P_{\ker (\g)} x P_{\ran (\e)^\perp}\eta_n}{\| \gamma_n P_{\ker (\g) \cap \ker (\guno)} \|     \| \eta_n\|}\\
  &    =  \frac{\gamma_n P_{\ker (\g) \cap \ker (\guno)} x P_{\ran (\e)}\eta_n}{\| \gamma_n P_{\ker (\g) \cap \ker (\guno)} \|     \| \eta_n\|}\\
  &    +\left (\frac{\| \gamma_n P_{\ker (\g)} \| \| P_{\ran (\e)^\perp} \eta_n\|}{ \| \gamma_n P_{\ker (\g) \cap \ker (\guno)}\| \| \eta_n \|} \right) \frac{ \gamma_n P_{\ker (\g)} x P_{\ran (\e)^\perp} \eta_n}{\| \gamma_n P_{\ker (\g)} \| \| P_{\ran (\e)^\perp} \eta_n\|}.
   \end{split}
  \end{equation*}
  Hence, depending on the value of the limit $\lambda$ of the sequence $\frac{\| \gamma_n P_{\ker (\g)} \| \| P_{\ran (\e)^\perp} \eta_n\|}{\| \gamma_n P_{\ker (\g) \cap \ker (\guno)}\| \| \eta_n\|}$, we have
       \begin{equation*}
\lim_n \mathbb{P} (\gamma_n x \eta_n) = \begin{cases} \mathbb{P}(\gdue x \e) & \mbox{ for } \lambda =0\\
\mathbb{P}( \guno x \euno) & \mbox{ for } \lambda = \infty\\
\mathbb{P} ( \gdue x \e +\lambda \guno x \euno) & \mbox{ for } \lambda \neq 0, \infty 
\end{cases}.
\end{equation*}
  
  Suppose now that we have
  \begin{equation*}
  \g x \edue =0, \quad \gdue x \e =0.
  \end{equation*}
  In this case, if $x$ has rank $2$, we must have $\ran (\e) = \ker (x)$ and $\ran (x) =\ker (\g)$ and so $y \in C_{\ker (\g), \ran (\e)}$.
  
Now consider the case
 \begin{equation*}
  \begin{array}{cc}
  \g x \e =0, & \g x \euno =0,\\
  \guno x \e =0, & \guno x \euno \neq 0,\\
  \g x \edue \neq 0, & \gdue x \e \neq 0. \end{array}
  \end{equation*}
  We shall see that this is impossible, since it would entail $\rank (x) =3$. For let $v_0 \in \ran (\e)$, $v_1 \in \ran (\euno)$ and $v_2 \in \ran (\edue)$ be non-zero vectors and note that $\rank (\gamma^{(i)}) = \rank (\eta^{(i)}) =1$ for every $i =0,1,2$ (since both $\gdue$ and $\edue$ are non-zero). From the hypothesis we should have $x v_0 \neq 0$, $x v_1 \neq 0$ and $x v_2 \neq 0$. But $x v_0 \in \ker (\g) \cap \ker (\guno)$, while $x v_1 \in \ker (\g) \backslash \ker (\guno)$; moreover $x v_2 \notin \ker (\g)$. Hence the range of $x$ contains three linearly independent vectors and its rank should be $3$.
  
  We are left to consider the case
   \begin{equation*}
  \begin{array}{cc}
  \g x \e =0, & \g x \euno =0,\\
  \guno x \e =0, & \guno x \euno \neq 0,\\
  \g x \edue \neq 0, & \gdue x \e = 0. \end{array}
  \end{equation*}
  In this instance $\ker (x) = \ran (\e)$ and we can write
  \begin{equation*}
  \begin{split}
  \frac{\gamma_n x \eta_n}{\| \gamma_n \| \| P_{\ran (\e)^\perp \cap \ran (\euno)^\perp} \eta_n\|} &=  \frac{\gamma_n x P_{\ran (\e)^\perp}  \eta_n}{\| \gamma_n \| \| P_{\ran (\e)^\perp \cap \ran (\euno)^\perp} \eta_n\|} \\
 &  =   \frac{\gamma_n x P_{\ran (\e)^\perp \cap \ran (\euno)^\perp} \eta_n}{\| \gamma_n \| \| P_{\ran (\e)^\perp \cap \ran (\euno)^\perp} \eta_n\|}
 +   \frac{\gamma_n x P_{\ran (\euno)} P_{\ran (\e)^\perp} \eta_n}{\| \gamma_n \| \| P_{\ran (\e)^\perp \cap \ran (\euno)^\perp} \eta_n\|}\\
 & =  \frac{\gamma_n x P_{\ran (\e)^\perp \cap \ran (\euno)^\perp} \eta_n}{\| \gamma_n \| \| P_{\ran (\e)^\perp \cap \ran (\euno)^\perp} \eta_n\|}\\
 &+
 \left( \frac{ \| P_{\ran (\e)^\perp} \eta_n \| \| \gamma_n P_{\ker (\g)} \|}{ \| P_{\ran (\e)^\perp \cap \ran (\euno)^\perp} \eta_n\| \|\gamma_n \|} \right) \frac{ \gamma_n P_{\ker (\g)} x P_{\ran (\euno) }P_{\ran (\e)^\perp} \eta_n}{\| \gamma_n P_{\ker (\g)} \|  \| P_{\ran (\e)^\perp} \eta_n \|}.
 \end{split}
  \end{equation*}
  Denote by $\lambda$ the limit of the sequence $\frac{\| P_{\ran (\e)^\perp} \eta_n \| \| \gamma_n P_{\ker (\g)} \| }{ \| P_{\ran (\e)^\perp \cap \ran (\euno)^\perp} \eta_n \| \| \gamma_n \|}$. We have
        \begin{equation*}
\lim_n \mathbb{P} (\gamma_n x \eta_n) = \begin{cases} \mathbb{P}(\g  x \edue) & \mbox{ for } \lambda =0\\
\mathbb{P}( \guno x \euno) & \mbox{ for } \lambda = \infty\\
\mathbb{P} (\g x \edue + \lambda \guno x \euno) & \mbox{ for } \lambda \neq 0, \infty 
\end{cases}.
\end{equation*}
  This situation is analogous to previous cases. To complete the proof note that convergence of the translates of $m_2$ is guaranteed, since the map from $M_{3,2}$ to $\mathbb{P}(M_3)$ is a homeomorphism with its image. $\Box$

     \begin{prop}
     \label{prop11}
     Let $\mu$ be a quasi-invariant top-dimensional finite Radon measure on $\partial_\beta G(\mathbb{Z})  $.Then the Koopman representation $\pi_\mu : C^*(G (\mathbb{Z}) \times G(\zz) ) \rightarrow  \mathbb{B}(L^2 (\mu))$ is a $c'_0 (G(\zz) \times G(\zz))$-representation.   
      \end{prop}
\proof Without loss of generality we can suppose that $\mu$ is a probability measure. In order to reach a contradiction, suppose that $\{g_n\} \subset \gz \times \gz$ is a diverging sequence and $\epsilon >0$, $\xi, \zeta \in L^2_1 (\partial_\beta \gz   , \mu)$ are given such that $  |\braket{g  \xi, \zeta}_\mu| >\epsilon$ for every $g \in \{g_n\}$, where we denote with $\braket{\cdot , \cdot}_\mu$ the scalar product associated to the measure $\mu$. We now apply the following reduction: by approximating $\xi$ and $\zeta$ with essentially bounded functions and changing the value of $\epsilon$ we can suppose that $\xi$ and $\zeta$ are in $L^\infty (\mu)$; now changing again the value of $\epsilon$ and using the inequality $|\braket{g  \xi, \zeta}_\mu| \leq \braket{g 1, 1}_\mu  \|\xi\|_\infty \|\zeta\|_\infty$ we can suppose that the diverging sequence $g_n=(g_n^{(1)}, g_n^{(2)}) \in \gz \times \gz$ satisfies $|\braket{g_n 1,1}_\mu  |>\epsilon$ for every $n \in \mathbb{N}$. \\
First we will work with the Radon Nykodin derivative of the push-forward of the measure $\mu$, considered as a measure on $\phi^{-1}(M_{3,2})$. Suppose that $g_n=(g_n^{(1)}, g_n^{(2)}) \in \gz \times \gz$ is such that $ |\braket{ g_n 1, 1 }_{\phi_* \mu}| > \epsilon$ for every $n \in \mathbb{N}$. If both components of $g_n$ diverge, up to taking a subsequence, by means of Lemma \ref{lemmapl}, we obtain that there are (line, plane)-sets $D_1^0,...,D_k^0$ and $D_1,...,D_l$ such that, denoting $D^0 = \bigcup_{i=1}^m D_i^0$ and $D=\bigcup_{j=1}^l D_j$, for every $m \in M_{3,2} \backslash D^0$, there exists the limit $\lim g_n^{(1)} m g_n^{(2)}$ and is an element of $D$ (whenever it belongs to $M_{3,2}$). 
  It follows from Lebesgue dominated convergence Theorem that 
 \begin{equation*}
  \lim_n g_n \phi_* \mu (f) = \lim_n \phi_* \mu (f \circ g_n)= 0
  \end{equation*}
   for every compactly supported continuous function $f$ which is zero on the finite union of linear sets $C =D \cup D^0$. In particular, for every compact set $K \subset M_{3,2} \backslash C$ we have $\lim_n g_n \phi_* \mu (\chi_K) =0$ and so, for every compact set $K' \subset \phi^{-1}(M_{3,2}) \backslash \phi^{-1} (C)$ we have $g_n \mu (\chi_{K'}) \leq g_n \mu (\chi_{\phi^{-1} (\phi (K'))}) = g_n \phi_* \mu (\phi (K')) \rightarrow 0$. Let $U \subset \phi^{-1}(M_{3,2})$ be an open set containing $\phi^{-1} (C)$ with $\mu (U) <\epsilon^2/9$ and $K \subset \phi^{-1}(M_{3,2}) \backslash U$ a compact set such that $\mu (\phi^{-1}(M_{3,2}) \backslash (U \cup K))< \epsilon^2/9$. We have, by an application of H{\" o}lder inequality,
 \begin{equation*}
 \begin{split}
 \braket{g_n 1, 1}_\mu &= \int_{\phi^{-1}(M_{3,2})} \sqrt{\frac{dg_n \mu}{d \mu}} d\mu \\
 &= \int_U \sqrt{\frac{dg_n \mu}{d \mu}} d\mu + \int_{\phi^{-1}(M_{3,2}) \backslash (U \cup K)}  \sqrt{\frac{dg_n \mu}{d \mu}}d \mu+ \int_K \sqrt{\frac{dg_n \mu}{d \mu}} d\mu\\
 &\leq \left(\int_U  \frac{dg_n \mu}{d \mu}d\mu\right)^{1/2} \mu (U)^{1/2} +  \left(\int_{\phi^{-1}(M_{3,2}) \backslash (U \cup K)}  \frac{dg_n \mu}{d \mu}d\mu\right)^{1/2} \mu  (\phi^{-1}(M_{3,2})\backslash (U \cup K))^{1/2} \\
 &+ \left(\int_K \frac{dg_n \mu}{d \mu} d\mu \right)^{1/2} \mu (K)^{1/2} \\
 & \leq \mu(U)^{1/2} + \mu (\phi^{-1}(M_{3,2})\backslash (U \cup K))^{1/2} + \left(\int_K \frac{dg_n \mu}{d \mu} d\mu \right)^{1/2} < 3 \sqrt{\epsilon^2/9}=\epsilon
 \end{split}
\end{equation*}
for $n$ large enough. $\Box$

 \subsection{Rational measures}
 \label{secrat}
 In this section we show that the Koopman representation associated to a rational measure on $\partial_\beta G(\mathbb{Z})$ is a $C^*_{red}$-representation. In order to do so, we will use a convenient description of $G(\mathbb{Z})$ as a dynamical system over the global stabilizer of a particular rational plane-line set. This dynamical interpretation reveals the amenability of the action of this global stabilizer on the boundary preimage of the corresponding plane-line set, entailing the amenability of $G(\mathbb{Z}) \times G(\mathbb{Z})$ on the rational part of the boundary.
 
 \begin{lem}
 \label{trans-rat}
The action of $G(\mathbb{Z})$ on the set of rational lines is transitive; the same is true for the set of rational planes.
\end{lem}
\proof Let $W$ be a rational line. Then there is a vector $(m,n,k) \in \mathbb{Z}^3 \cap V$. Suppose that $m$, $n$ and $k$ are all non-zero. Let $q$ be the greatest common divisor between $n$ and $k$ and write $n=q N$, $k =qK$. If follows from B\'{e}zout identity that there are $a$, $b \in \mathbb{Z}$ such that $Na-Kb=1$. Letting
\begin{equation*}
g=\left( \begin{array}{cc}	N	&	b	\\
				K	&	a	\end{array}\right) \quad \in \SL(2,\mathbb{Z}),
				\end{equation*}
				we have $g (1,0)^t = (N,K)^t$. Hence $g^{-1} (N,K)^t = (1,0)^t$ and so $g^{-1} (n,k)^t = (q,0)^t$. We thus have
 \begin{equation*}
\left( \begin{array}{cc}	1	&\begin{array}{cc}	0	&	0	\end{array}	\\
				\begin{array}{c}0	\\	0\end{array}	&	g^{-1}	\end{array}\right) \left(\begin{array}{c}	m	\\	n	\\	k	\end{array}\right) = \left( \begin{array}{c}	m	\\	q	\\	0\end{array}\right).
				\end{equation*}
 Proceeding in the same way, we find $h \in \SL(2,\mathbb{Z})$ and $u \in \mathbb{Z}$ such that $h (m,q)^t = (u, 0)$ and so
 \begin{equation*}
 \left( \begin{array}{cc}	h	& \begin{array}{c}0	\\	0\end{array}	\\
				\begin{array}{cc}	0	&	0	\end{array}	&	1	\end{array}\right) 
 \left( \begin{array}{cc}	1	&\begin{array}{cc}	0	&	0	\end{array}	\\
				\begin{array}{c}0	\\	0\end{array}	&	g^{-1}	\end{array}\right) \left(\begin{array}{c}	m	\\	n	\\	k	\end{array}\right) = \left( \begin{array}{c}	u	\\	0	\\	0\end{array}\right).
				\end{equation*}
Clearly, if one among $m$, $n$ and $k$ is zero, the same strategy applies.\\
Let now $W$ be a rational plane. Note that for every $g \in G(\mathbb{Z})$, we have $(gW)^\perp = {}^t{g}^{-1} (W^\perp)$. Hence transitivity of the action on rational planes is guaranteed by the first part of the proof. $\Box$

Let $R =\braket{e_1, e_2}$ and $R^\perp=\braket{e_3}$. We have
\begin{equation*}
H_\lambda :=\gz_R = \left\{ \left(\begin{array}{cc}	g	&	 \begin{array}{c}s	\\ t \end{array}\\
						\begin{array}{cc}	0	&	0	\end{array}	&	1	\end{array}\right) \; | \; g \in SL(2,\mathbb{Z}) ,\; s,t \in \mathbb{Z}\right\} ,
						\end{equation*}
						\begin{equation*}
H_\rho:= \gz_{R^\perp} =\left\{ \left(\begin{array}{cc}	f	&	 \begin{array}{c}0	\\ 0 \end{array}\\
						\begin{array}{cc}	h	&	k	\end{array}	&	1	\end{array}\right) \; | \; f \in SL(2,\mathbb{Z}) ,\; h,k \in \mathbb{Z}\right\} .
						\end{equation*}
In particular, $(\gz \times \gz )_{\phi^{-1} (C_{V,W})} = H_\lambda \times H_\rho$. Note that the set $C_{R, R^\perp}$ is the subset of $M_{3,2}$ of pairs of projective matrices $(m_1,m_2)$ with the property that $m_1$ has the form
\begin{equation*}
m_1 = \left( \begin{array}{ccc} \times	&	\times	&	0\\
		\times	&	\times	&	0\\
		0	&	0	&	0\end{array}\right) \qquad \mod \rr\backslash 0.
		\end{equation*}
We will make use of the map
 \begin{equation*}
 \wedge : \begin{array}{ccc}\mathbb{Z}^2 \times \mathbb{Z}^2 &\rightarrow& M_2 (\mathbb{Z}) \\
 	&	&	\\
 			\left( \left(\begin{array}{c}v_1\\v_2\end{array}\right) , \left(\begin{array}{c}w_1\\w_2\end{array}\right)\right) & \mapsto & v\wedge w=\left( \begin{array}{cc} v_1 w_1 	&	v_1 w_2\\ v_2 w_1 & v_2 w_2 \end{array}\right). \end{array}
			\end{equation*}
			\label{eqrel}
			This map satisfies the following properties for $g  \in \SL(2,\mathbb{Z})$, $u,v,z \in \zz^2$:
			\begin{equation}
				\label{eqrel}
				\begin{array}{cc}
			g(u \wedge v) = (gu) \wedge v,  & (u \wedge v) g = u \wedge ({}^t g v), \\ (u+v)\wedge z = u\wedge z + v \wedge z, & u \wedge (v+z)= u\wedge v + u \wedge z.\end{array}
			\end{equation}
The (plane, line)-set we are interested in is $C_{R,R^\perp}$. We shall now give a dynamical characterization of $\gz$ as $H_\lambda \times H_\rho$-space. For let $  0 \neq l \in \mathbb{Z}$ and consider the space
 \begin{equation*}
 X_l := \{ (x_l , m , n) \; | \; x_l \in \GL (2,\mathbb{Q}), \; \det (x_l)=l^{-1}, \; m,n \in \mathbb{Z}^2, \;x_l + l^{-1} m \wedge n \in M_2 (\mathbb{Z})\}.
 \end{equation*}
 We want to define an action of $H_\lambda \times H_\rho$ on $X_l$; it is convenient to use the following compact form for representing elements in the acting group:  
 \begin{equation*}
 \mbox{if }\quad \gamma=\left(\begin{array}{cc}	g	&	 \begin{array}{c}s	\\ t \end{array}\\
						\begin{array}{cc}	0	&	0	\end{array}	&	1	\end{array}\right) \in H_\lambda,  \quad v=\left(\begin{array}{c} s\\ t \end{array}\right) \qquad \mbox{ we write }\quad \gamma= g \ltimes v;
						\end{equation*}
\begin{equation*}
 \mbox{if } \quad \eta =  \left(\begin{array}{cc}	f	&	 \begin{array}{c}0	\\ 0 \end{array}\\
						\begin{array}{cc}	h	&	k	\end{array}	&	1	\end{array}\right), \quad w= \left(\begin{array}{c} h\\k\end{array}\right) \qquad \mbox{ we write }\quad \eta = f   \bar{\ltimes}  w. \end{equation*}
  
 Let $H_\lambda$ act on $\mathbb{Z}^2$ through $(g \ltimes v)_l(u) := g u + lv$ and similarly let $H_\rho$ act on $\mathbb{Z}^2$ through $(f \bar{\ltimes} w)_l (u) = {}^t f^{-1} (u-lw) $. The action of $H_\lambda \times H_\rho$ on $X_l$ is given by
 \begin{equation*}
( (g \ltimes v), (f \bar{\ltimes} w)) \cdot  (x_l, m , n) := (g x_l f^{-1}, (g\ltimes v)_l (m), (f \bar{\ltimes} w)_l (n)).
\end{equation*}
For $l=0$ we let
\begin{equation*}
X_0 := \{ (x_{i,j})_{i,j=1}^3 \in \gz \; | \; x_{3,3}=0\}.
\end{equation*}
 
 \begin{lem}
 The $H_\lambda \times H_\rho$-space $\gz$ is equivariantly isomorphic to $\bigsqcup_{l \in \mathbb{Z}} X_l$.
 \end{lem}
 \proof
 Let $y \in \gz$ and write
 \begin{equation*}
	y= \left(\begin{array}{cc}	x	&	m	\\	
						{}^tn	&	l	\end{array}\right), \quad \mbox{ with } x \in M_2 (\mathbb{Z}), \; m \in \mathbb{Z}^2, \; n \in \mathbb{Z}^2, l \in \mathbb{Z}.
						\end{equation*}
Note that $X_0$ is $H_\lambda \times H_\rho$-invariant.  \\
Consider the map
  \begin{equation*}
  \alpha  : \begin{array}{ccc}\gz \backslash X_0 &\rightarrow &\bigsqcup_{l \in \mathbb{Z} \backslash \{0\}} X_l \\
  	&	&	\\
 				\left(\begin{array}{cc}	x	&	m	\\	
						{}^tn	&	l	\end{array}\right) & \mapsto &  (x - l^{-1} m \wedge n , m , n).
\end{array} \end{equation*}
A direct computation, using the relations (\ref{eqrel}), shows that $\alpha$ is $H_\lambda \times H_\rho$-equivariant. By definition, $\alpha$ is injective. In order to check that this map actually takes values in $\bigsqcup_{l \in \mathbb{Z}} X_l$, for
\begin{equation*}
x=\left(\begin{array}{cc}	x_{1,1}	&	x_{1,2}	\\
				x_{2,1}	&	x_{2,2}	\end{array}\right),
				\end{equation*}
 write
\begin{equation*}
\tilde{x} = \left(\begin{array}{cc}	x_{2,2}	&	-x_{1,2}	\\
				-x_{2,1}	&	x_{1,1}	\end{array}\right);
\end{equation*}	
then we have, using the formula $\det(x + x') = \det (x) + \det (x') + \Tr (\tilde{x} x')$ for $x$, $x' \in M_2 (\mathbb{R})$,			
\begin{equation*}
\det (x-l^{-1} m \wedge n) = \det(x) - l^{-1} \Tr ([\tilde{x} (m)] \wedge n) = \det (x) - l^{-1} \braket{\tilde{x} (m), n},
\end{equation*}
 while
 \begin{equation*}
 \det(y)=l \det(x)- \Tr ([\tilde{x}(m)] \wedge n)= l \det (x) - \braket{\tilde{x} (m), n } =1.
 \end{equation*}
 Hence
 \begin{equation*}
 \det (x-l^{-1} m \wedge n) = \det(x)-l^{-1} (l \det(x) -1)=l^{-1}.
 \end{equation*}
 In order to show that $\alpha$ is surjective, let $l \neq 0$ and $(x_l , m , n) \in X_l$, then
 \begin{equation*}
 \begin{split}
 \det \left( \begin{array}{cc}	x_l +l^{-1} m\wedge n	&	m	\\
 			{}^t n 	&	l	\end{array}\right) &= l \det(x_l) +\det (x_l) \Tr([x_l^{-1}(m)] \wedge n) - \det(x_l) \Tr([x_l^{-1}(m)] \wedge n) \\ &= l \det(x_l) =1.\end{split}
			\end{equation*} 
 This completes the proof. $\Box$\\
 
 
In order to prove amenability of the action of $H_\lambda \times H_\rho$ on $\phi^{-1} (C_{V,W})$, we proceed in the following way: first we show amenability of the action on $\bigsqcup_l \partial X_l$; secondly we prove that the part of $\phi^{-1} (C_{R,R^\perp})$ which does not lie in  $\bigsqcup_l \partial X_l$ is an amenable $H_\lambda\times H_\rho$-space.
 
 \begin{lem}
 \label{lemz}
Let $l \in \mathbb{Z} \backslash \{0\}$. There exists a sequence $\mu^{(\lambda,l)}_n : \mathbb{Z}^2 \rightarrow \mathcal{P}(H_\lambda)$ with the following property: for every $\gamma=g\ltimes v \in H_\lambda$, $\epsilon >0$ there are $N \in \mathbb{N}$ and $K \subset \mathbb{R}^2$ compact such that
\begin{equation*}
\| \gamma\mu^{(\lambda,l)}_N (u) - \mu^{(\lambda,l)}_N (\gamma_l u)\|_1 < \epsilon \qquad \mbox{ for every }  u \notin lK \cap \mathbb{Z}^2.
\end{equation*}
The compact set $K$ does not depend on the choice of $l$. A similar result holds for $H_\rho$.
\end{lem}
 \proof Consider the following maps:
 \begin{itemize}
 \item[(i)] $\alpha : \mathbb{Z}^2   \rightarrow \mathbb{R}^2 $, $\alpha (v)=v/l$. This map is $H_\lambda$-equivariant if we consider on $\mathbb{Z}^2$ the $l$-th action and on $\mathbb{R}^2$ the canonical action of $H_\lambda$ given by $(g \ltimes v)(u )= gu + v$.
 \item[(ii)] $\beta : \mathbb{R}^2 \backslash \{0\} \rightarrow \mathbb{R}P^1$, $\beta (m, n) = m/n$. This map is $\SL(2,\mathbb{Z})$-equivariant. Fix a metric $d$ on $\mathbb{R}P^1$ which induces the natural topology, for every $g \ltimes v \in H_\lambda$ we have $\lim_{x_n \rightarrow \infty}d ( \beta((g \ltimes v)(x_n)), \beta (g(x_n)))= 0$ (c.f. \cite{ozae}, 2).
 \item[(iii)] An approximate invariant mean $\gamma_n : \mathbb{R}P^1 \rightarrow \mathcal{P}(\SL(2,\mathbb{Z}))$ (indexed by $\mathbb{N}$) given by amenability of the action of $\SL(2,\mathbb{Z})$ on $\mathbb{R}P^1$.
 \item[(iv)] An approximate invariant mean $\delta_\lambda : \mathcal{P} (\Delta_\beta \SL(2,\mathbb{Z})) \rightarrow \mathcal{P}(H_\lambda)$ (indexed by a net $\Lambda$) given by the amenability of the action of $H_\lambda$ on $\mathcal{P}(\Delta_\beta \SL(2,\mathbb{Z}))$. Here the action on $\Delta_\beta \SL(2,\mathbb{Z})$ is intended to be the extension of the natural action coming from the identification $\SL(2,\mathbb{Z}) = H_\lambda /\mathbb{Z}^2$. Note that $g \ltimes v (h) = gh$ for every $g \ltimes v \in H_\lambda$, $h \in \SL(2,\mathbb{Z})$.
  \end{itemize}
 Let $\epsilon >0$, $F \subset H_\lambda$ finite and $u \in \mathbb{Z}^2 \backslash (\{0\} \cup \bigcup_{(g \ltimes v) \in F} \{-lg^{-1} (v)\})$. For $\lambda \in \Lambda$, $n \in \mathbb{N}$ we have  \begin{align*}
  \| (g  \ltimes v)& \delta_\lambda \gamma_n \beta \alpha (u) - \delta_\lambda \gamma_n \beta \alpha ((g \ltimes v)_l (u))\|_1\\
  &\leq \| (g \ltimes v) \delta_\lambda \gamma_n \beta \alpha (u) - \delta_\lambda (g \gamma_n \beta \alpha (u))\|_1 & (1)  \\
  &+ \| \delta_\lambda (g \gamma_n \beta \alpha (u)) - \delta_\lambda \gamma_n (g \beta \alpha (u))\|_1 & (2)   \\
  &+ \| \delta_\lambda \gamma_n (g \beta \alpha (u)) - \delta_\lambda \gamma_n \beta \alpha (gu)\|_1  & (3)  \\
 &+ \| \delta_\lambda \gamma_n \beta \alpha (gu) - \delta_\lambda \gamma_n \beta \alpha ((g \ltimes v)_l(u))\|_1   & (4) \end{align*}  
There is $\bar{\lambda} \in \Lambda$ such that contribution (1) is less than $\epsilon /3$ for $\lambda \geq \bar{\lambda}$, uniformly in $u$ and $n$, for $g  \ltimes v \in F$. Fix such a $\lambda$. There is $\bar{n} \in \mathbb{N}$ such that contribution (2) is smaller than $\epsilon /3$ for every $n >\bar{n}$ (by continuity of $\delta_\lambda$) uniformly on $u$, for $g  \ltimes v \in F$. Fix such an $n$. Contribution (3) is identically $0$ by $\SL(2,\mathbb{Z})$-equivariance of $\alpha$ and $\beta$. There is a compact set $K \subset \mathbb{R}^2$ (which we assume to contain $\{0\} \cup \bigcup_{(g \ltimes v) \in F} \{-g^{-1} (v)\}$) such that contribution (4) is smaller than $\epsilon /3$ for every $\alpha (u) \notin K$ (by continuity of $\delta_\lambda \gamma_n \beta$), hence for every $u \notin l K$ and $g  \ltimes v \in F$. We set $\mu^{(\lambda,l)}_\epsilon = \delta_\lambda \gamma_n \beta \alpha $ on $\mathbb{Z}^2 \backslash \{0\}$ and extend arbitrarily on $\{0\}$.\\Taking an exhaustion of $H_\lambda$ by finite sets $F_n\subset F_{n+1}$ and a sequence $\epsilon_n \rightarrow 0 $, we obtain a sequence of compact sets $K_n \subset K_{n+1} \subset \mathbb{R}^2$ and of maps $\mu^{(\lambda,l)}_n:=\mu^{(\lambda,l)}_{\epsilon_n} : \mathbb{Z}^2 \rightarrow \mathcal{P}(H_\lambda)$ with the property that
\begin{equation*}
\| \mu_n ((g \ltimes v)_l(u)) - (g\ltimes v)(\mu_n (u)) \|_1 < \epsilon_n \qquad \mbox{ for every } u \notin lK_n, \; g \in F_n.
\end{equation*}
The first part of the claim is proved. In order to show the same result for $H_\rho$, note that the group isomorphism $H_\rho \rightarrow H_\lambda$, $f \bar{\ltimes} w \mapsto {}^t (f \bar{\ltimes} w)^{-1}$ interchanges the $H_\rho$ and $H_\lambda$ actions on $\mathbb{Z}^2$. $\Box$

 For the proof of the next result we need to recall a construction given in \cite{Sk}. Let $\mathbb{H}$ be the hyperbolic plane endowed with its hyperbolic metric and for $a$, $b \in \mathbb{H}$ let 
 \begin{equation*}
 F_{a,b} = \{ x \in \mathbb{H} \; | \; d (x  , [a,b])<1\}.
 \end{equation*}
 Fix a reference point $x_0 \in \mathbb{H}$. For $g \in \SL(2,\mathbb{R})$ let $\chi' (g)$ be the characteristic function of $F_{x_0 , g(x_0)}$ and $\chi (g) = \chi'(g) / \| \chi' (g)\|_2 \in L^2 (\mathbb{H})$. Arguing as in \cite{Sk}, proof of Teorem 4.4, we note that for every $\epsilon >0$ and $F \subset \SL(2,\mathbb{R}) \times \SL(2,\mathbb{R})$ finite, there is a compact set $K \subset \SL(2,\mathbb{R})$ such that
 \begin{equation}
 \label{eqeq3}
 \| \chi_{g h f^{-1}} - \gamma \chi_h \|_2 < \frac{\epsilon}{2} \qquad \mbox{ for every } h \notin K, \; (g, f) \in F.
 \end{equation}

 \begin{prop}
 \label{propbamen1}
 For every $l\in \mathbb{Z}$ the action of $H_\lambda \times H_\rho$ on $\partial_\beta X_l$ is amenable.
 \end{prop}
 \proof First we consider the case $l=0$. Let $x \in X_0$ and write
 \begin{equation*}
 x=\left( \begin{array}{cc} x_0 	& h\\
 		{}^t k 	& 	0 \end{array} \right), \qquad \mbox{ with } h,k \in \zz^2.
		\end{equation*}
Let $ g\ltimes v \in H_\lambda$, $f \bar{\ltimes} w \in H_\rho$ be such that $((g\ltimes v), (f \bar{\ltimes} w)) x=x$, then $g h=h$ and ${}^t f^{-1} k =k$, which shows that the stabilizer of $x$ is an amenable group. Thus the action of $H_\lambda \times H_\rho$ on $X_0$ is boundary amenable.
 
 Fix $l \in \mathbb{Z}$, $l \neq 0$. In virtue of Lemma \ref{lemz}, given an exhaustion of $H_\lambda $ by finite \textit{symmetric} sets $F_n \subset F_{n+1}$, we can find finite subsets $C_n \subset C_{n+1}\subset \mathbb{Z}^2$ and maps $\mu_n : \mathbb{Z}^2 \rightarrow \mathcal{P}(H_\lambda )$ such that
 \begin{equation*}
 \| \mu_n ( (g \ltimes v)_l (u)) - (g\ltimes v) \mu_n (u)\|_1 < 1/n^2 \qquad \mbox{ for every } u \notin C_n, \; (g \ltimes v) \in F_n.
 \end{equation*}
  We can suppose that $F_n C_n \subset C_{n+1}$ for every $n$. Let now $D_n := F_n C_n$; then the above estimate is still true after replacing $C_n$ with $D_n$. The point is that the sets $D_n \backslash D_{n-1}$ are "almost" eventually $F$-invariant for every finite subset $F\subset H_\lambda$, in the sense that for every finite set $F$ and $n$ large enough we have $F (D_n \backslash D_{n-1}) \subset D_{n+1} \backslash D_{n-2}$. Consider the map $h: \mathbb{Z}^2 \rightarrow \mathbb{N}$ given by $h(u):=\min\{n \in \mathbb{N} \; | \; u \in D_n\}$ and let $\mu : \mathbb{Z}^2  \backslash D_1 \rightarrow \mathcal{P}(H_\lambda)$ be given by (c.f. \cite{ozae}, 3)
  \begin{equation*}
  \mu (u) = \frac{1}{h(u)} \sum_{i=1}^{h(u)-1} \mu_i (u).
  \end{equation*}
 Extend arbitrarily $\mu$ on the whole $\mathbb{Z}^2$. For every $F \subset H_\lambda$ finite and $\epsilon >0$ there is a finite set $K \subset \mathbb{Z}^2$ such that
 \begin{equation}
 \label{eqba}
 \| \mu ((g \ltimes v)_l (u))- (g \ltimes v) \mu (u)\|_1 < \epsilon \qquad \mbox{ for every } u \notin K, \; (g\ltimes v) \in F.
 \end{equation}
Relation (\ref{eqba}) is exactly the boundary amenability of the $l$-th action of $H_\lambda$ on $\mathbb{Z}^2$ (\cite{ozak}, Proposition 4.1). Hence there is an isometry $U: l^2 (\mathbb{Z}^2) \rightarrow l^2 (H_\lambda)$ such that
\begin{equation*}
U^* \lambda (s) U - \alpha_{\zz^2}(s) \in \mathbb{K} (l^2 (\mathbb{Z}^2)) \qquad \mbox{ for all } s \in H_\lambda,
\end{equation*} 
where we denote by $\alpha_{\zz^2}$ the unitary representation of $H_\lambda$ on $l^2(\zz^2)$.
  Consider now the countable discrete space $\tilde{X}_l$ given by triples $(x_l, g\ltimes v, u)$ with $x_l \in \GL(2,\mathbb{Q})$ such that $l x_l \in \GL(2,\mathbb{Z})$, $\det(x_l) =1/l$, $g\ltimes v \in H_\lambda$, $u \in \mathbb{Z}^2$ endowed with the action of $H_\lambda \times H_\rho$ given by
 \begin{equation*}
 (g \ltimes v, f \bar{\ltimes} w) (x_l , g' \ltimes v' , u) = (g x_l f^{-1}, (g \ltimes v)( g' \ltimes v'), ( f \bar{\ltimes} w)_l (u)).
 \end{equation*}
 Since this action is proper, it is implemented by a representation $\pi_0$ of the reduced group $C^*$-algebra of $H_\lambda \times H_\rho$; in particular it is boundary amenable. It follows from \cite{ozak} Proposition 4.1 that there is an isometry $W : l^2(\tilde{X}_l) \rightarrow l^2  (H_\lambda \times H_\rho)$ such that
 \begin{equation}
 \label{eqeq2}
 W^* \lambda (s)W - \pi_0 (s) \in \mathbb{K} (l^2 (\tilde{X}_l)) \qquad \mbox{ for every } s \in H_\lambda \times H_\rho.
 \end{equation}
Let $V: l^2 (X_l ) \rightarrow l^2 (\tilde{X}_l)$ be given by $V(\delta_{x_l} \otimes \delta_m \otimes \delta_n) = \delta(x_l) \otimes U\delta_m \otimes \delta_n$. Arguing as in \cite{ozak} Proposition 4.1 ($3 \Rightarrow 2$), we have 
\begin{equation}
\label{eqeq1}
\lim_{m \rightarrow \infty}\|  \pi_0 (s) V (\delta_{x_l} \otimes \delta_m \otimes \delta_n) - V \alpha_{ {X}_l} (s )(\delta_{x_l} \otimes \delta_m \otimes \delta_n)\|_2 =0  \qquad \mbox{ for all }   s \in H_\lambda \times H_\rho
\end{equation} 
 Let then $\phi^{(\lambda)} : X_l \rightarrow l^2_1 (H_\lambda \times H_\rho)$, $\phi^{(\lambda)} ((x_l,m,n)) = W V(\delta_{x_l} \otimes \delta_m \otimes \delta_n)$. Let $s \in H_\lambda \times H_\rho$ and $x=(x_l , m,n) \in X_l$; then
 \begin{equation*}
 \begin{split}
 \|\lambda(s)  \phi^{(\lambda)} (x) -  \phi^{(\lambda)} (sx)\|_2 & = \|\lambda(s) WV\delta_x - WV \delta_{sx}\|_2 \\
 	& \leq \| \lambda(s) WV (\delta_{x_l} \otimes \delta_m \otimes \delta_n) - W \pi_0 (s) V (\delta_{x_l} \otimes \delta_m \otimes \delta_n) \|_2 \\ &+ \| W \pi_0 (s) V (\delta_{x_l} \otimes \delta_m \otimes \delta_n) - WV \alpha_{X_l} (\delta_{x_l} \otimes \delta_m \otimes \delta_n)\|_2 .
	\end{split}
   \end{equation*}
 Since $V$ is an isometry, it follows from (\ref{eqeq2})  and (\ref{eqeq1}) that 
 \begin{equation*}
 \lim_{m \rightarrow \infty}   \|\lambda(s)  \phi^{(\lambda)} ((x_l,m,n)) -  \phi^{(\lambda)} (s(x_l,m,n))\|_2 =0 \qquad \mbox{ for every } s \in H_\lambda \times H_\rho.
 \end{equation*}
 Putting $ \mu^{(\lambda)} := |  \phi^{(\lambda)}|^2 : X_l \rightarrow l^1 (H_\lambda \times H_\rho)$, we have
 \begin{equation}
 \label{eqeq22}
 \lim_{m \rightarrow \infty}   \| s  \mu^{(\lambda)} ((x_l,m,n)) -  \mu^{(\lambda)} (s(x_l,m,n))\|_1 =0 \qquad \mbox{ for every } s \in H_\lambda.
 \end{equation}

In a similar way, we find $ \mu^{(\rho)} : X_l \rightarrow l^1 (H_\lambda \times H_\rho)$ such that
\begin{equation}
\label{eqeq7}
\lim_{n \rightarrow \infty} \|s  \mu^{(\rho)} ((x_l,m,n)) -  \mu^{(\rho)} (s(x_l,m,n))\|_1 =0 \qquad \mbox{ for every } s \in H_\lambda \times H_\rho.
\end{equation}

Let now $\chi : \SL(2,\mathbb{R}) \rightarrow L^2 (\mathbb{H})$  be the map introduced above. Consider the isometry $V'' : l^2 (X_l) \rightarrow L^2 (\mathbb{H}) \otimes l^2 (X_l)$, $\delta_{(x_l , m, n )} \mapsto \chi_{  \sign(l)\sqrt{|l|}\cdot x_l  } \otimes \delta_{(x_l , m ,n)}$. Let $\mathbb{H} \times X_l$ be endowed with the proper action of $H_\lambda \times H_\rho$ given by
\begin{equation*}
((g \ltimes v), (f \bar{\ltimes} w)) (y, x_l , m,n) = (g y , g x_l f^{-1} , (g \ltimes v)_l (m), (f \bar{\ltimes} w)_l (n)).
\end{equation*}
 By properness, the Koopman representation on $L^2 (\mathbb{H}) \otimes l^2 (X_l)$ is a representation $\pi_0$ of $C^*_{red} (H_\lambda \times H_\rho)$. In virtue of (\ref{eqeq3}), for every $\epsilon >0$ and $F \subset H_\lambda \times H_\rho$ finite there is a compact set $K \subset \SL(2,\mathbb{R})$ such that
 \begin{equation}
 \label{eqeq4}
 \| ((V'')^* \pi_0 (s) V'' -\alpha_{X_l} (s)) (\delta_{(x_l , m ,n)}) \|_2 <\frac{\epsilon}{4} \qquad \mbox{ for every } x_l \notin K/\sign(l) \sqrt{|l|}, \; s \in F.
 \end{equation}
 By Voiculescu's Theorem there is an isometry $U'' : L^2 (\mathbb{H}) \otimes l^2 (X_l) \rightarrow l^2 (H_\lambda \times H_\rho)$ such that
 \begin{equation*}
 (U'')^* \lambda  (s) U'' - \pi_0 (s) \in \mathbb{K} (L^2 (\mathbb{H}) \otimes l^2 (X_l)) \qquad \mbox{ for every } s \in H_\lambda \times H_\rho.
 \end{equation*}
Hence, for every $\epsilon >0$ and finite set $F \subset H_\lambda \times H_\rho$ there is a finite set $K' \subset X_l$ such that
\begin{equation*}
\| ((U'')^* \lambda  (s) U'' - \pi_0 (s)) (\chi_{\sign(l) \sqrt{|l|} \cdot x_l} \otimes \delta_{(x_l,m,n)})\|_2 < \frac{\epsilon}{4} \qquad \mbox{ for every } (x_l,m,n) \notin K', \; s \in F.
\end{equation*}  
 Projecting on the first component, we find that for every $\epsilon >0$ and finite set $F \subset H_\lambda \times H_\rho$ there is a finite set $K_1 \subset M_2(\mathbb{Z}/l)$ such that
\begin{equation}
\label{eqeq5}
\| ((U'')^* \lambda  (s) U'' - \pi_0 (s)) (\chi_{\sign(l) \sqrt{|l|} \cdot x_l} \otimes \delta_{(x_l,m,n)})\|_2 < \frac{\epsilon}{4} \qquad \mbox{ for every } x_l \notin K_1, \; s \in F.
\end{equation}
Let $\epsilon >0$ and $F \subset H_\lambda \times H_\rho$. Pick compact sets $K \subset \SL(2,\mathbb{R})$ as in (\ref{eqeq4}) and $K_1 \subset M_2 (\mathbb{Z}/l)$ as in (\ref{eqeq5}). Then, for $x_l \notin K/\sign(l) \sqrt{|l|} \cup K_1$ and $s \in F$ we have
\begin{equation*}
\begin{split}
\|\lambda & (s) U'' V'' (\delta_{(x_l,m,n)}) - U'' V'' \alpha_{X_l}(s) (\delta_{(x_l,m,n)}) \|_2\\
&\leq \| \lambda  (s) U'' V'' (\delta_{(x_l,m,n)}) - U'' \pi_0 (s) V'' (\delta_{(x_l,m,n)}) \|_2\\& + \| U'' \pi_0 (s) V'' (\delta_{(x_l,m,n)}) - U'' V'' \alpha_{X_l}(s) (\delta_{(x_l,m,n)})\|_2\\
&\leq  \| \lambda  (s) U'' \chi_{  \sign(l)\sqrt{|l|} \cdot x_l  } \otimes \delta_{(x_l , m ,n)} - U'' \pi_0 (s) \chi_{  \sign(l)\sqrt{|l|}  \cdot x_l } \otimes \delta_{(x_l , m ,n)} \|_2 \\ &+ \| \pi_0 (s) V'' (\delta_{(x_l,m,n)}) -   V'' \alpha_{X_l} (\delta_{(x_l,m,n)})\|_2 < \epsilon / 2 + \epsilon /2 = \epsilon.
\end{split}
\end{equation*}
It follows that the map $\mu^{(X_l)} : X_l \rightarrow l^1 (H_\lambda \times H_\rho)$, $(x_l , m, n) \mapsto |U'' V''  \delta_{(x_l,m,n)}|^2 $ satisfies:  
\begin{equation}
\label{eqeq6}
\lim_{x_l \rightarrow \infty} \| s \mu^{(X_l)} ((x_l,m,n)) - \mu^{(X_l)} (s (x_l,m,n)) \|_1 =0 \qquad \mbox{ for every } s \in H_\lambda \times H_\rho.
\end{equation}
 
 Let now $F_n \subset F_{n+1}$ be an exhaustion of $\SL(2,\mathbb{Z} ) \times \SL(2,\mathbb{Z})$ by finite symmetric sets and let $C_n$ be an exhaustion of $M_2 (\mathbb{Z}/l)$ by finite sets such that for every $n \in \mathbb{N}$ we have $F_n C_n \subset C_{n+1}$ and define $D_n = F_n C_n$; let $h^{(X_l)} : M_2 (\mathbb{Z}/l) \rightarrow \mathbb{N}$ be given by $h^{(X_l)}(x_l) = \min\{n \in \nn \; | \; x_l \in D_n\}$. Similarly define exhaustions of $H_\lambda$, $H_\rho$, $\mathbb{Z}^2$ and functions $h^{(\lambda)}$, $h^{(\rho)}$. Define $\phi : X_l \rightarrow l^1 (H_\lambda \times H_\rho)$, $\mu ((x_l,m,n)) = (h^{(X_l)} (x_l) + h^{(\lambda)} (m)+h^{(\rho)} (n) )^{-1} (h^{(X_l)}(x_l) \mu^{(X_l)} ((x_l,m,n)) + h^{(\lambda)}(m)  \mu^{(\lambda)} ((x_l,m,n)) +h^{(\rho)} (n)  \mu^{(\rho)} ((x_l,m,n)))$. We have
 \begin{equation*}
 \lim_{(x_l,m,n) \rightarrow \infty} \|s\mu ((x_l,m,n)) - \mu(s(x_l,m,n)) \|_1 =0 \qquad \mbox{ for every } s \in H_\lambda \times H_\rho.
 \end{equation*}
 This proves that the action of $H_\lambda \times H_\rho$ on $\partial_\beta (X_l)$ is amenable. $\Box$

 \begin{lem}
 \label{lem9}
 For every $l \in \mathbb{Z}$ and $k \in \mathbb{N}$ let $X_{l,k}:= \{(x_l, m ,n) \in X_l \; | \; \min\{\|m\|, \|n\|\} \leq k|l|\}$. We have $\phi^{-1} (C_{R, R^{\perp}}) \cap( \partial_\beta \gz \backslash \bigsqcup_{l \in \mathbb{Z}} \partial_\beta X_l) \subseteq \partial_\beta \gz \backslash \bigcup_{k \in \mathbb{N}} \partial_\beta (\bigsqcup_{l \in \mathbb{Z}} X_{l,k})$.
 \end{lem}
 \proof
 Observe the following general fact: if $A$ is an infinite subset of $\gz$, $\omega \in \partial_\beta A \subset \gz$ and $\{x_\lambda \} \subset \gz$ is a net converging to $\omega$, then $x_\lambda$ eventually belongs to $A$; in order to see this, let $x_\lambda \rightarrow \omega$, this means that the principal ultrafilters associated to $x_\lambda$ (the principal ultrafilter associated to an element $x$ is the ultrafilter given by all the subsets containing $x$) converge to $\omega$; since $\omega \in \partial_\beta A$, by definition, $\omega$ belongs to the open set determined by $A$, i.e. the set of ultrafilters containing $A$; then, eventually, the principal ultrafilters associated to $x_\lambda$ belong to this open set and so $A$ is eventually an element of these ultrafilters; it follows that the points $x_\lambda$ eventually belong to $A$. On the other hand, if $\omega \in \partial_\beta  \gz\backslash \partial_\beta A$, then $x_\lambda$ eventually belongs to $A^c$, since $ \gz\backslash \partial_\beta A= \{ \omega \in \partial_\beta \gz \; | \; A^c \in \omega\}$.\\
  It follows that if $\omega \in  \partial_\beta \gz \backslash \bigsqcup_{l \in \mathbb{Z}} \partial_\beta X_l$, any net $x_\lambda  $ converging to $\omega$ should satisfy $|(x_\lambda)_{3,3}| \rightarrow \infty$.\\
 Let now $x_n$ be a sequence in $\gz$ such that $|(x_n)_{3,3}| \rightarrow \infty$ and $\phi (x_n) \rightarrow \bar{x} \in C_{R,R^{\perp}}$. Let $x$ be a representative of $\bar{x}$ in $M_{3, \|\cdot \|=1}$ (abusing notation, we are replacing $\bar {x} = (m_1,m_2)$ by its first component $m_1$) and $p$ be the orthogonal projection on $R$, then we have $p x p = x \neq 0$. It follows from Lemma \ref{lemams} that there is $\theta \in \mathbb{R} \backslash \{0\}$ such that $p x_n p/\|p x_n p\| \rightarrow \theta x$. We write
 \begin{equation*} 
 x_n = \left( \begin{array}{cc}	g_n	&	v_n\\
 				{}^t w_n	&	l_n	\end{array}\right), \qquad x=\left( \begin{array}{cc}	g	&\begin{array}{c}0\\ 0 \end{array}	\\
				\begin{array}{cc} 0 & 0\end{array}  & 0\end{array}\right),
				\end{equation*}
 where $v_n$, $w_n \in \mathbb{Z}^2$. Note that since $x$ has rank $2$, then $g $ is invertible. Let $\lambda_n = \| p x_n p\| \theta$, so that $x_n/ \lambda_n \rightarrow x$, which entails $g_n/ \lambda_n \rightarrow g$, $v_n /\lambda_n \rightarrow 0$, $w_n /\lambda_n \rightarrow 0$ and $l_n /\lambda_n \rightarrow 0$. Denote
 \begin{equation*}
 \tilde{g}_n = \left( \begin{array}{cc}	d_n	&-b_n \\
 				-c_n	&	a_n	\end{array}\right), \qquad \mbox{ where } \quad g_n = \left(\begin{array}{cc}	a_n	&	b_n	\\
				c_n	&	d_n \end{array}\right).
				\end{equation*}
 We can write $ 1=\det (x_n) = l_n \det (g_n) - \braket{ \tilde{g}_n (v_n), w_n}$ for every $n$; hence $\det (g_n) = (1+ \braket{ \tilde{g}_n (v_n), w_n})/l_n$. Since both $\lambda_n$ and $|l_n|$ go to infinity, we have $(\lambda^2_n l_n)^{-1} \rightarrow 0$; moreover $\tilde{g}_n /\lambda_n$ converges since the same is true for $g_n /\lambda_n$. But then, the relation
 \begin{equation*}
 0 \neq |\det (g)| = \lim_n \left| \det (g_n/\lambda_n)\right| = \lim_n\left| \braket{ \frac{\tilde{g}_n}{\lambda_n} \left(\frac{v_n}{ |l_n \lambda_n|^{1/2}}\right), \frac{w_n}{|l_n \lambda_n|^{1/2}}}\right|
 \end{equation*} 
 implies
 \begin{equation}
 \label{eq1}
 \frac{v_n}{|l_n \lambda_n|^{1/2}} \rightarrow v_0 \neq 0, \qquad  \frac{w_n}{|l_n \lambda_n|^{1/2}} \rightarrow w_0 \neq 0.
 \end{equation}
 We will see how equation (\ref{eq1}) entails the claim. Let $\omega \in \phi^{-1} (C_{R, R^{\perp}}) \cap( \partial_\beta \gz \backslash \bigsqcup_{l \in \mathbb{Z}} \partial_\beta X_l)$ and $x_\lambda \rightarrow \omega$. We already observed that $|(x_\lambda)_{3,3}| =|l_\lambda|\rightarrow \infty$. Suppose for the sake of contradiction that $\omega \in  \bigcup_{k \in \mathbb{N}} \partial_\beta (\bigsqcup_{l \in \mathbb{Z}} X_{l,k})$. By the first part of the proof there should be $k \in \mathbb{N}$ such that $x_\lambda$ eventually belongs to $\bigsqcup_{l \in \mathbb{Z}} X_{l,k}$. Since $M_{3,2}$ is first countable, we can extract a sequence $x_n$ from the net $x_\lambda$ such that $\phi (x_n) \rightarrow \phi (\omega)$ and $|l_n| \rightarrow \infty$. But then the relations $|v_n| \leq |l_n| k$ and $l_n/\lambda_n \rightarrow 0$ would entail $v_n / |l_n \lambda_n|^{1/2} \rightarrow 0$ and similarly for $w_n$, contradicting equation (\ref{eq1}). $\Box$

 \begin{lem}
 \label{lem10}
 Let $\Gamma$ be a countable discrete exact group and $X$ a countable discrete $\Gamma$-space. Let $Y \subset \partial_\beta X$ be a closed $\Gamma$-invariant subset and $\eta : X \rightarrow \mathcal{P}(  \Gamma)$ be such that for every $\gamma  \in \Gamma$, $\omega \in Y$ and net $x_\lambda \rightarrow \omega$ we have
 \begin{equation*}
 \lim_\lambda  \|\eta (\gamma x_\lambda) -  \gamma \eta (x_\lambda)\|_1 =0 .
 \end{equation*}
 Then the action of $\Gamma$ on $Y$ is amenable.
 \end{lem}
  
 \proof (c.f. \cite{ozak} Proposition 4.1) Let $\tilde{\eta} : \Delta_\beta X \rightarrow \mathcal{P}(\Delta_\beta \Gamma)$ be the continuous extension of $\eta$ to $\Delta_\beta X$ ($\mathcal{P}(\Delta_\beta  \Gamma)$ is endowed with the $w^*$-topology on the dual of $l^\infty ( \Gamma)$). Since $ \Gamma$ is exact, the action of $ \Gamma$ on $\Delta_\beta  \Gamma$ is amenable. Hence it is enough to show that $\tilde{\eta}|_Y$ is equivariant. For let $\omega \in Y$, $\gamma \in  \Gamma$ and $\{x_\lambda\}_{\lambda \in \Lambda}$ be a net in $X$ such that $x_\lambda \rightarrow \omega$. Let $\epsilon >0$ and $f \in l^\infty ( \Gamma)$; there is $\bar{\lambda} \in \Lambda$ such that
 \begin{equation*}
 | \gamma \tilde{\eta}(\omega)(f) - \gamma \eta (x_\lambda)(f)| < \epsilon /3, \qquad | \eta (\gamma x_\lambda)(f) -   \tilde{\eta}(\gamma \omega)(f)| < \epsilon /3,
 \end{equation*}
 \begin{equation*}
 \|\eta (\gamma x_\lambda) -  \gamma \eta (x_\lambda)\|_1 < \epsilon /(3 \|f\|) \qquad \mbox{ for } \lambda \geq \bar{\lambda},
 \end{equation*}
 where the first two inequalities follow from the continuity of $\eta$ and the third one follows from the assumption on $Y$ in the statement. Hence, for $\lambda \geq \bar{\lambda}$, we have
 \begin{equation*}
 \begin{split}
 | \gamma \eta(\omega) (f) - \eta(\gamma \omega)(f)| \leq& | \gamma \eta(\omega)(f) - \gamma \eta (x_\lambda)(f)| + | \gamma \eta (x_\lambda)(f) - \eta (\gamma x_\lambda)(f)|\\ & + | \eta (\gamma x_\lambda)(f) - \eta (\gamma \omega)(f)| < \epsilon.
 \end{split}
 \end{equation*}
 Since $\epsilon$ and $f$ are arbitrary, the result follows. $\Box$

 \begin{prop}
 \label{propbamen2}
 The action of $H_\lambda \times H_\rho$ on $\partial_\beta \gz \backslash (\bigcup_{k \in \mathbb{N}} (\partial_\beta \bigsqcup_{l \in \mathbb{Z}} X_{l,k}) \cup (\bigsqcup_{l \in \mathbb{Z}} \partial_\beta X_l))$ is amenable.
 \end{prop}
 \proof 

 It follows from Lemma \ref{lemz} that there are an exhaustion of $H_\lambda \times H_\rho$ by finite set $F_n \subset F_{n+1}$, a sequence of natural numbers $k_n \rightarrow \infty$ and sequences of maps $\mu_n^{(\lambda,l)} : \mathbb{Z}^2 \rightarrow \mathcal{P}(H_\lambda)$, $\mu_n^{(\rho,l)} : \mathbb{Z}^2 \rightarrow \mathcal{P}(H_\rho)$ such that for every $n \in \mathbb{N}$ and $(\gamma_1, \gamma_2) \in F_n$ we have
 \begin{equation*}
 \|\gamma_1 \mu^{(\lambda,l)}_n (u) - \mu^{(\lambda,l)}_n ((\gamma_1)_l u)\|_1 < \frac{1}{2n^2} \qquad \mbox{ for every } \|u\| > |l|k_n, \; l \in \mathbb{Z} \backslash 0,
 \end{equation*}
 \begin{equation*}
 \|\gamma_2 \mu^{(\rho,l)}_n (v) - \mu^{(\rho,l)}_n ((\gamma_2)_l v)\|_1 < \frac{1}{2n^2} \qquad \mbox{ for every } \|v\| > |l| k_n, \; l \in \mathbb{Z}\backslash 0.
 \end{equation*} 
Take the sequences of map $\mu^{(\lambda)}_n$, $\mu^{(\rho)}_n$ from $G(\zz) = \bigsqcup_l X_l$ to $\mathcal{P}(H_\lambda)$ and $\mathcal{P}(H_\rho)$ respectively given by $\mu^{(\lambda)}_n ((x_l,u,v)) = \mu^{(\lambda,l)}_n (u)$ and $\mu^{(\rho)}_n ((x_l,u,v)) = \mu^{(\rho,l)}_n (v)$ for $l \neq 0$, arbitrary for $l=0$.
 
We define a sequence of maps $\eta_n : \gz \rightarrow \mathcal{P}(H_\lambda \times H_\rho)$ by
 \begin{equation*}
  \eta_n   = \mu^{(\lambda)}_n  \otimes \mu_n^{(\rho)} .
 \end{equation*}
We have
 \begin{equation}
 \label{eqesaust}
 \begin{split}
\| s\eta_n (x)& - \eta_n (sx)\|_1 \leq \|\gamma_1 \mu^{(\lambda)}_n (u) - \mu^{(\lambda)}_n ((\gamma_1)_l u)\|_1 +  
  \|\gamma_2 \mu^{(\rho)}_n (v) - \mu^{(\rho)}_n ((\gamma_2)_l v)\|_1 \\
  &< \frac{1}{n^2}\qquad \mbox{ for every } x \notin \bigsqcup_{l \in \mathbb{Z}\backslash \{0\}}X_{l,k_n}, \; s=(\gamma_1,\gamma_2) \in F_n.
\end{split}
\end{equation}
 Up to taking a subsequence we can suppose that $k_{n+1} > k_n^2$ for every $n$. Let $s \in H_\lambda \times H_\rho$. We will see now that for every $\gamma$ in either $H_\lambda$ or $H_\rho$ there is a positive constant $M=M(\gamma)$ such that for every $u \in \mathbb{Z}^2$ we have:
 \begin{equation}
 \label{stimam}
 \frac{1}{M}\| u \| - |l|M \leq \| \gamma_l u\| \leq M \|u\| +|l|M.
 \end{equation}
 In order to check this, note that it is enough to show the existence of positive real numbers $m_1,m_2,m_3,$ and $m_4$ such that for every $u \in \mathbb{R}^2$ we have
 \begin{equation*}
 m_1 \|u\| - m_2|l| \leq \|\gamma_lu\| \leq m_3 \|u\| + m_4|l|;
 \end{equation*}
suppose then that $\gamma = g \ltimes v \in H_\lambda$; in this case we can take $m_2=m_4 = \|v\|$, $m_3 = \|g\|$. It is left to show the existence of $m_1 >0$ such that $\| gu\| \geq m_1 \|u\|$ for every $u \in \mathbb{R}^2 \backslash \{0\}$. Suppose such an $m_1$ does not exist, then there is a sequence $u_n$ such that $\|gu_n\|/\|u_n\| \rightarrow 0$. Up to taking a subsequence we can suppose that $u_n /\|u_n\|$ converges to a point $u$ of norm one; hence $g u = \lim_n g (u_n /\|u_n\|) =0$, which is impossible. Suppose now that $\gamma = f \bar{\ltimes} w \in H_\rho$. In this case we can take $m_2=m_4= \| {}^t f^{-1}\| \|w\|$, $m_3= \|{}^t f^{-1}\|$ and for $m_1$ we can apply the above reasoning to ${}^t f^{-1}$ in place of $g$. This proves the desired relation.\\
We want to show that as a consequence of (\ref{stimam}), since we are assuming $k_{n+1}^2 > k_n$, there is $N=N(s) \in \mathbb{N}$ such that
 \begin{equation}
  \label{eq10}
  s(  \bigsqcup_{l\neq 0} X_{l,k_{n+1}}\backslash \bigsqcup_{l \neq 0} X_{l,k_n}) \subset \bigsqcup_{l \neq 0} X_{l,k_{n+2} } \backslash \bigsqcup_{l \neq 0} X_{l,k_{n-1} } \qquad \mbox{ for every } n>N.
  \end{equation}
  In fact, let $s = (\gamma_1 , \gamma_2)$ and suppose $u \in \mathbb{Z}^2$ is such that $|l| k_n \leq \|u\| \leq |l| k_{n+1}$ for some $n$. Letting $M=M(\gamma_1)$, in virtue of (\ref{stimam}), we have
  \begin{equation*}
  \|(\gamma_1)_l u\| \leq M \|u\| + |l| M \leq |l|M (k_{n+1} +1) \leq |l| k_{n+1}^2 < |l| k_{n+2} \qquad \mbox{ for } 2M < k_{n+1},
  \end{equation*}
  \begin{equation*}
  \|(\gamma_1)_lu\| \geq \frac{1}{M} \| u\| - M |l| \geq \frac{1}{M} |l| k_n - M |l| \geq |l| \sqrt{k_{n}} > |l|k_{n-1} \qquad \mbox{ for } 2M < \sqrt{k_n}.
  \end{equation*}
  Hence, for $n$ large enough we will have $|l|k_{n-1} < \|(\gamma_1)_l u \| < |l| k_{n+2}$. In the same way we obtain $|l|k_{n-1} < \|(\gamma_2)_l u \| < |l| k_{n+2}$ for $n$ large enough. This proves (\ref{eq10}).\\
 For every $n \in \mathbb{N}$ define $\eta|_{\bigsqcup_{l \in \mathbb{Z}}  X_{l,k_{n+1}} \backslash  \bigsqcup_{l \in \mathbb{Z}} X_{l,k_n}}: \bigsqcup_{l \in \mathbb{Z}}  X_{l,k_{n+1}} \backslash  \bigsqcup_{l \in \mathbb{Z}} X_{l,k_n} \rightarrow \mathcal{P}(H_\lambda \times H_\rho)$ by
\begin{equation*}
\eta (x) = \frac{1}{n} \sum_{i=1}^{n} \eta_i (x) 
\end{equation*}
and extend to $\gz$, arbitrarily inside $\bigsqcup_{l \in \mathbb{Z}}  X_{l,k_1}$. Now $s \in F_m$ for some $m \in \mathbb{N}$ (hence $s$ belongs to $F_n$ for every $n>m$) and so equation \ref{eqesaust} is true for $s$ for every $n\geq m$. 
In virtue of equation (\ref{eq10}), for every $x \in   \bigsqcup_{l \neq 0}  X_{l,k_{n+1}} \backslash  \bigsqcup_{l \neq 0} X_{l,k_n} $, with $n >\max\{N,m\}$, we have the following possibilities, depending on where $x$ lands under the action of $g$:
 \begin{equation*}
 \| g\eta (x) - \eta (gx) \|_1 = \| \frac{1}{n}\sum_{i=1}^{n} \eta_i (x) - \frac{1}{n} \sum_{i=1}^{n} \eta_i (gx)\|_1 \leq \frac{2N +1}{n},
 \end{equation*}
 \begin{equation*}
  \| g\eta (x) - \eta (gx) \|_1 = \| \frac{1}{n} g\sum_{i=1}^n \eta_i (x) - \frac{1}{n-1} \sum_{i=1}^{n-1} \eta_i (gx)\|_1 \leq \frac{N+1}{n  } +\frac{N}{n-1}+\frac{1}{n} ,
  \end{equation*}
  \begin{equation*}
 \| g\eta (x) - \eta (gx) \|_1 = \| \frac{1}{n} g\sum_{i=1}^n \eta_i (x) - \frac{1}{n+1} \sum_{i=1}^{n+1} \eta_i (gx)\|_1 \leq \frac{N}{n}+\frac{N+1}{n+1}+\frac{1}{n+1}.
 \end{equation*}
Hence for every $\epsilon >0$ there is $N' (s,\epsilon) \in \mathbb{N}$ such that 
\begin{equation}
\label{eqsepsilon}
 \| g\eta (x) - \eta (gx) \|_1 < \epsilon \qquad \mbox{ for every } x \notin \bigsqcup_{l \neq 0} X_{l,k_{N'(s,\epsilon)}}.
 \end{equation}
  Let now $\omega \in  \partial_\beta \gz \backslash (\bigcup_{k \in \mathbb{N}} (\partial_\beta \bigsqcup_{l \in \mathbb{Z}} X_{l,k}) \cup (\bigsqcup_{l \in \mathbb{Z}} \partial_\beta X_l))$ and $x_\lambda \rightarrow \omega$. As we observed in the proof of Lemma \ref{lem9} we must have $(x_\lambda)_{3,3} \neq 0$ eventually and $x_\lambda$ should eventually avoid every set of the form $\bigsqcup_{l \in \mathbb{Z}} X_{l,k}$. So for every $n \in \mathbb{N}$ there is $\bar{\lambda}$ in the net such that $x_\lambda \notin \bigsqcup_{l\neq 0} X_{l,k_n}$ for every $\lambda \geq \bar{\lambda}$. Hence for every $\epsilon >0$ there is $\bar{\lambda}$ such that $x_\lambda \notin \bigsqcup_{l\neq 0} X_{l,k_{N' (s,\epsilon)}}$ for every $\lambda \geq \bar{\lambda}$ and so equation \ref{eqsepsilon} entails
  \begin{equation*}
  \| s \eta (x_\lambda) - \eta (s x_\lambda)\|_1 < \epsilon \qquad \mbox{ for every } \lambda \geq \bar{\lambda};
  \end{equation*}
  since $\epsilon$ is arbitrary we have $\lim_\lambda \| g \eta (x_\lambda) - \eta (g x_\lambda)\|_1 =0$. The result follows from Lemma \ref{lem10}. $\Box$

 \begin{lem}
 \label{lemamen}
 Let $\Gamma$ be a countable discrete exact group acting on a compact space $X$. Let $X = Y \sqcup Z$ be a decomposition of $X$ in $\Gamma$-invariant sets, with $Y$ closed and $Z$ open. Suppose that the restricted actions of $\Gamma$ on $Y$ and $Z$ are both amenable. Then the action of $\Gamma$ on $X$ is amenable.
 \end{lem}
 \proof By exactness, $C(X) \rtimes_r \Gamma / C_0 (Z) \rtimes_r \Gamma = C(Y) \rtimes_r \Gamma$. Since $ C_0 (Z) \rtimes_r \Gamma$ and $C(Y) \rtimes_r \Gamma$ are both nuclear, also $C(X) \rtimes_r \Gamma$ is nuclear. Hence the action of $\Gamma$ on $X$ is amenable in virtue of \cite{An1} Theorem 3.4. $\Box$

  \begin{prop}
  \label{rat}
 Let $\mu$ be a $G(\mathbb{Z}) \times G(\mathbb{Z})$-quasi-invariant finite rational measure on $\partial_\beta G(\zz)$. Then $\pi_\mu$ is a $C^*_{red}(G(\zz) \times
 G(\zz))$-representation.
 \end{prop}
 \proof In virtue of Proposition \ref{propbamen1} the action of $H_\lambda \times H_\rho$ on the open set $\bigsqcup_l X_l$ is amenable; in virtue of Lemma \ref{lem9} and Proposition \ref{propbamen2} the action of $H_\lambda \times H_\rho$ on the closed set $\phi^{-1} (C_{R,R^\perp}) \cap (\partial_\beta G(\zz) \backslash \bigsqcup_l \partial_\beta X_l)$ is amenable; hence, by Lemma \ref{lemamen} the action of $H_\lambda \times H_\rho$ on $\phi^{-1} (C_{R, R^\perp}) \cap \partial_\beta G(\zz)$ is amenable. It follows from Lemma \ref{lemamencoset} and Lemma \ref{trans-rat} that the action of $G(\zz) \times G(\zz)$ on the inverse image under $\phi$ of the disjoint union of all the rational (plane, line)-sets is amenable. It follows from Lemma \ref{wk} that the Koopman representation associated to $\mu$ factors through $C^*_{red} (G(\zz) \times G(\zz))$. $\Box$

 \subsection{Irrational measures}
 \label{secirr}
 Parallel to the previous section, in the following we want to prove amenability of the action of $G(\zz) \times G(\zz)$ on the orbits of preimages of irrational (plane, line)-sets. Actually, in this case we obtain a stronger result, namely the global stabilizers of irrational (plane, line)-sets act amenably on $\partial_\beta G(\zz)$ (and not just on the preimages of the corresponding (plane, line)-sets). 
  \begin{prop} 
 \label{irrational}
 Let $V$ be either an irrational line or an irrational plane. Then $G (\mathbb{Z})_V$ is amenable.
 \end{prop}
 \proof 
Let $V$ be an irrational line and $v \in V$ non-zero. Consider the $\mathbb{Q}$-algebra $A_v := \{ a \in M_3 (\mathbb{Q}) \; | \; a V \subset V\}$. For every $a \in A_v$ there is a real number $\lambda (a)$ such that $av=\lambda (a) v$. The map $A_v \rightarrow \mathbb{R}$, $a\mapsto \lambda(a)$ is a character. Up to permutation of the indices, we can suppose $v=(v_1,v_2,v_3)^t$ with $v_1$ and $v_2$ independent over $\mathbb{Q}$. \\Consider first the case in which $v_1, v_2$ and $v_3$ are linearly independent over $\mathbb{Q}$; then $\lambda$ is injective and the global stabilizer $\gz_V$ is abelian.\\
Consider now the case in which $v_3= \theta_1 v_1 + \theta_2 v_2$ is a $\mathbb{Q}$-linear combination of $v_1$ and $v_2$. Every element $a$ in the kernel of $\lambda$ must be of the form
\begin{equation*}
a=\left( \begin{array}{ccc}	\alpha_1 \theta_1 	&	\alpha_1 \theta_2 	& -\alpha_1	\\
				\alpha_2 \theta_1 	& 	\alpha_2 \theta_2 	& 	-\alpha_2	\\
				\alpha_3 \theta_1 	& 	\alpha_3 \theta_2 	&-\alpha_3	\end{array}\right), \qquad \mbox{ for some } \alpha_1, \alpha_2, \alpha_3 \in \mathbb{Q}.
				\end{equation*} 
In particular, for such $a$'s, we have $a (e_1+ \theta_1 e_3)=  a(e_2 + \theta_2 e_3)=0$; on the other hand any matrix $a$ satisfying these relations must satisfy $\lambda (a)=0$. Denote by $p_2$ the orthogonal projection onto the (real) linear span of $e_1 + \theta_1 e_3$ and $e_2 + \theta_2 e_3$. Hence
\begin{equation*}
\ker (\lambda) =\{ a \in A_v \; | \; a p_2 =0\}.
\end{equation*}
It follows that $\ker (\lambda) \cap p_2 A_v p_2 =\{0\}$ and so the restriction of $\lambda$ to $p_2 A_v p_2$ is injective, entailing commutativity of $p_2 A_v p_2$. Note now that, since $\ker (\lambda)$ is a bilateral ideal and $(1-p_2) \in \ker (\lambda)$, for every $a \in A_v$ we have $(1-p_2) a p_2=0$. Up to a change of basis, performed by an orthogonal  matrix in $\SL(3,\mathbb{Q})$, we can suppose that $p_2$ projects onto the span of $e_1$ and $e_2$; the relation $(1-p_2) a p_2=0$ translates then to the fact that $a$ has the form
\begin{equation}
\label{eqq}
\left( \begin{array}{ccc}	\times	&	\times	&	\times	\\
			\times	&	\times	&	\times	\\
			0	&	0	&	\times	\end{array}\right).
			\end{equation}
If $a \in M_3 (\mathbb{Q})$ is invertible and has the form (\ref{eqq}), then also $p_2 a p_2$ is invertible. Hence the group of invertibles in $A_v$ is a subgroup of the semidirect product group of the form $\Inv (p_2 A_v p_2) \ltimes H$, where $H$ is the group of matrices of the form
\begin{equation*}
\left(\begin{array}{ccc}	1	&	0	&	x\\
					0	&	1	&	y\\ 
					0	&	0	&	z\end{array}\right),\qquad x,y,z \in \mathbb{Q}, \; z \neq 0.
					\end{equation*}
Since $H$ and $\Inv (p_2 A_v p_2)$ are amenable groups, it follows that the group of invertibles in $A_v$ is amenable as well.\\ 
The case in which $V$ is a rational plane follows as well, indeed in this instance the global stabilizer of $V $ is the transposed of the global stabilizer of $V^\perp$, which is amenable. $\Box$

\begin{lem}
\label{exa}
Let $G$ be a discrete countable group and $H \subset G$ an amenable subgroup. If $G$ is exact, the left-right action of $G \times H$ and $H \times G$ on $\partial_\beta G$ are amenable.
\end{lem}
\proof 
In virtue of \cite{ozak} Prop. 4.1 (3) it is enough to prove that the action of $G \times H$ on $G$ is amenable. For let $\mu_n^H : \{*\} \rightarrow \mathcal{P} (H)$ and $\mu_n^G : G/H \rightarrow \mathcal{P}(G)$ be approximate invariant means. Then the maps $x \mapsto \mu_n^G (xH) \otimes \mu_n^H$, $G \rightarrow \mathcal{P}(G \times H)$ constitute an approximate invariant mean. The other case is similiar. $\Box$

\begin{prop}
\label{irr}
Let $\mu$ be a quasi-invariant irrational measure on $\partial_\beta G(\zz)$. Then the Koopman representation $\pi_\mu : C^*(G(\mathbb{Z}) \times G(\mathbb{Z})) \rightarrow \mathbb{B}(L^2 (\mu))$ factors through $C^*_{red} (G(\mathbb{Z}) \times G(\mathbb{Z}))$. 
\end{prop}
\proof The measure $\mu$ is supported on a union of $G(\zz) \times G(\zz)$-orbits of (plane, line)-sets. In virtue of Lemma \ref{ideal} we need to show the result for the orbit of an arbitrary irrational (plane, line)-set $C_{V,W}$. It follows from Proposition \ref{irrational} and Lemma \ref{exa} that the global stabilizer of $C_{V,W}$ acts amenably on $\partial_\beta \gz$. It follows from Lemma \ref{lemamencoset} and Lemma \ref{wk} that the Koopman representation $\pi_\mu$ is a $C^*_{red} (G(\zz) \times G(\zz))$-representation. $\Box$

\subsection{From (plane, line)-sets to (line, plane)-sets}
\label{s5}
We summarize the results of sections \ref{secrat} and \ref{secirr} in the following
\begin{prop}
\label{atom}
Let $\mu$ be a $G(\mathbb{Z}) \times G(\mathbb{Z})$-quasi-invariant finite Radon measure on $\partial_\beta G(\mathbb{Z})$ supported on the orbit of the preimage of a (plane, line)-set. Then the Koopman representation $\pi_\mu : C^* (G(\mathbb{Z}) \times G(\mathbb{Z})) \rightarrow \mathbb{B}(L^2 (\mu))$ is a $C^*_{red} (G(\mathbb{Z}) \times G(\mathbb{Z}))$-representation.
\end{prop}
\proof Follows from Proposition \ref{rat} and Proposition \ref{irr}. $\Box$

Now we want to see how a similar result holds when considering measures supported on orbits of preimages of (line, plane)-sets.

 \begin{lem}
 \label{lem4.18}
 Let $A$ be a set in either $\mathcal{E}$ or $\mathcal{F}$ and $\mu$ be a $G(\mathbb{Z}) \times G(\mathbb{Z})$-quasi-invariant Radon measure supported on the orbit of $\phi^{-1} (A)$ which annihilates every set in $\phi^{-1}(\mathcal{C})$. Then the Koopman representation associated to $\mu$ is a $C^*_{red} (G(\mathbb{Z}) \times G(\mathbb{Z}))$-representation.
 \end{lem}
     \proof We let $\{e_1, e_2, e_3\}$ be the standard orthonormal basis for $\mathbb{R}^3$. It follows from Lemma \ref{globalstabsets} that the global stabilizer of $A$ is of the form $G(\mathbb{Z})_V \times G(\mathbb{Z})_{V'}$, where either $\dim (V)=\dim (V')=1$ or $\dim(V)=\dim (V')=2$. Hence, if at least one between $V$ and $V'$ is irrational, the global stabilizer of $A$ acts amenably on $\phi^{-1} (A)$ in virtue of Proposition \ref{irrational} and Lemma \ref{exa}. Moreover, if both $V$ and $V'$ are rational, by transitivity (Lemma \ref{trans-rat}), the orbit of $\phi^{-1}(A)$ under $G(\mathbb{Z}) \times G(\mathbb{Z})$ is the orbit of the preimage of either the set of matrices whose range is $\braket{ e_1,e_2}$ and kernel contained in $\braket{e_1, e_2}$ in the case $\dim (V)=\dim (V')=2$, or the set of matrices whose range contains $e_3$ and kernel is $\mathbb{R} e_3$. The global stabilizer of $F_{\braket{ e_1,e_2}, \braket{e_1,e_2}}$ is $H_\lambda \times H_\lambda$ and the global stabilizer of $E_{\mathbb{R} e_3, \mathbb{R} e_3}$ is $H_\rho \times H_\rho$.\\
We want to show that both the action of $H_\lambda \times H_\lambda$ and of $H_\rho \times H_\rho$ on $\partial_\beta G(\mathbb{Z})$ is amenable. Consider first the case of $H_\lambda \times H_\lambda$. Let $G(\mathbb{Z}) = H_\lambda \sqcup (G(\mathbb{Z}) \backslash H_\lambda)$. Then the action of $H_\lambda \times H_\lambda$ on $\partial_\beta H_\lambda \subset \partial_\beta G(\mathbb{Z})$ is amenable by \cite{ozae}. Let now $x \in G(\mathbb{Z}) \backslash H_\lambda$ and we write
\begin{equation*}
x= \left( \begin{array}{cc} x_0 & m\\
			n	&	l\end{array}\right), \qquad \mbox{ where } x_0 \in M_2 (\mathbb{Z}), \; m,n \in \mathbb{Z}^2, \; l \in \mathbb{Z}.
			\end{equation*}   
 If $g \ltimes v $ and $ h \ltimes w$ in $H_\lambda$ are such that $(g \ltimes v) x (h \ltimes w)=x$, then $h^t (n) =n$, hence the stabilizer of $x$ is a subgroup of the product of $G(\mathbb{Z})$ and an amenable group; since $G(\mathbb{Z}) \backslash H_\lambda$ is invariant under taking inverses, it follows that the same reasoning can be applied to the equality $(h \ltimes w)^{-1} x^{-1} (g \ltimes v)^{-1} =x^{-1}$ in order to show that also $h \ltimes v$ belongs to an amenable group (fixed by x) and so the stabilizer of every point in $G(\mathbb{Z}) \backslash H_\lambda$ is amenable. It follows that the action of $H_\lambda \times H_\lambda$ on $G(\mathbb{Z}) \backslash H_\lambda$ is amenable and so the action of $H_\lambda \times H_\lambda$ on $G(\mathbb{Z})$ is boundary amenable. \\
 We shall see now that the boundary amenability of the action of $H_\rho \times H_\rho$ on $G(\mathbb{Z})$ follows from that of $H_\lambda \times H_\lambda$. In fact, let $\mu : G(\mathbb{Z}) \rightarrow \mathcal{P}(H_\lambda \times H_\lambda)$ be such that for every $g=(g_1,g_2) \in H_\lambda \times H_\lambda$ we have $\lim_{x \rightarrow \infty} \| g\mu (x) - \mu(g_1 x g_2^{-1})\|_1 =0$ and consider the map $\tilde{\mu} : G(\mathbb{Z}) \rightarrow \mathcal{P}(H_\rho \times H_\rho)$ given by $\tilde{\mu} (x) (\gamma) =  \mu ({}^t x^{-1})({}^t \gamma^{-1})$. Then for every $g=(g_1,g_2) \in H_\rho \times H_\rho$ we have $\lim_{x \rightarrow \infty} \| g \tilde{\mu} (x) - \tilde{\mu} (g_1 x g_2^{-1}) \|_1 =0$. The result follows from \cite{ozak} Proposition 4.1.\\
Let $X$ be either $F_{\braket{ e_1,e_2}, \braket{e_1,e_2}}$ or $E_{\mathbb{R} e_3, \mathbb{R} e_3}$ and $\Lambda$ its global stabilizer (hence either $H_\lambda \times H_\lambda$ or $H_\rho \times H_\rho$), let $\{g_i\}_{g_i \in G(\mathbb{Z} \times G(\mathbb{Z}))}$ be a set of coset representatives for $(G(\mathbb{Z}) \times G(\mathbb{Z}))/\Lambda$ and consider the space $(G(\mathbb{Z}) \times G(\mathbb{Z}))/\Lambda \times \phi^{-1} (X)$ endowed with the action of $G(\mathbb{Z}) \times G(\mathbb{Z})$ given by $g (g_i \Lambda, x)= (g_k \Lambda, \eta x)$, where $g g_i \Lambda = g_k \Lambda$ and $\eta = g_k^{-1} g g_i$. This dynamical system satisfies the conditions of Lemma \ref{lemamencoset} (with $Y=(\Lambda, \phi^{-1} (X))$) and so the action of $G(\mathbb{Z}) \times G(\mathbb{Z})$ on this space is amenable. The space $(G(\mathbb{Z}) \times G(\mathbb{Z})) \phi^{-1} (X)$ is a $G(\mathbb{Z}) \times G(\mathbb{Z})$-equivariant quotient of this space, the quotient map being given by $(g_i \Lambda, x) \mapsto g_i x$. Let now $\mu$ be a measure as in the hypothesis. Since the intersections of (line, line)-sets and (plane-plane)-sets are (plane, line)-sets (see the discussion at the beginning of Section \ref{section4}), it follows from the assumptions that the Koopman representation associated to the pullback measure on $(G(\mathbb{Z}) \times G(\mathbb{Z}))/\Lambda \times \phi^{-1} (X)$ and the one associated to $\mu$ are unitarily equivalent. It follows that $\pi_\mu$ is a $C^*_{red} (G(\mathbb{Z}) \times G(\mathbb{Z}))$-representation. $\Box$

 \begin{lem}
 \label{lem4.19}
 Let $D$ be a set in $\mathcal{D}$ and $\mu$ be a $G(\mathbb{Z}) \times G(\mathbb{Z})$-quasi-invariant Radon measure supported on the orbit of $\phi^{-1} (D)$ which annihilates every set in $\phi^{-1} (\mathcal{E} \cup \mathcal{F})$. Then the Koopman representation associated to $\mu$ is a $C^*_{red} (G(\mathbb{Z}) \times G(\mathbb{Z}))$-representation. 
 \end{lem}
 \proof We proceed as in the proof of Lemma \ref{lem4.18}. Hence we only need to show that the global stabilizer $\Lambda$ of a set $D$ in $\mathcal{D}$ acts amenably on its preimage in $\partial_\beta G(\mathbb{Z})$. In order to do so it is enough to show that for every $V, W \subset \mathbb{R}^3$ with $\dim(V)=1$, $\dim(W)=2$, we have that $G(\mathbb{Z})_V \times G(\mathbb{Z})_W$ acts amenably on $\pi(D_{V,W})$, where $\pi$ is the map sending a matrix to the cartesian product of its range and its kernel (this space is endowed with the quotient topology coming from $\SL(3,\mathbb{R}) \times \SL(3,\mathbb{R})$). In order to do so, observe that the action of $G(\mathbb{\mathbb{R}})_V \times G(\mathbb{\mathbb{R}})_W$ is transitive on $D_{V,W}$ (since the action of $\SL(3,\mathbb{R})$ is transitive on the Grassmanian); in particular it is amenable, since the stabilizer of $(\braket{e_1 , e_2}, \mathbb{R} e_3)$ in $\pi(D_{\mathbb{R} e_1, \braket{e_2, e_3}})$ is amenable, being given by products of Borel subgroups. The result follows. $\Box$
 
 We summarize the results of this section in the following
 \begin{prop}
 \label{proposition4.20}
 Let $\mu$ be a $G(\mathbb{Z}) \times G(\mathbb{Z})$-quasi-invariant finite Radon measure on $\partial_\beta G(\mathbb{Z})$ supported on the orbit of the preimage of a (line, plane)-set. Then the Koopman representation $\pi_\mu : C^* (G(\mathbb{Z}) \times G(\mathbb{Z})) \rightarrow \mathbb{B}(L^2 (\mu))$ is a $C^*_{red} (G(\mathbb{Z}) \times G(\mathbb{Z}))$-representation.
\end{prop}
\proof The measure $\mu$ can be decomposed in the orthogonal sum of a measure supported on the orbit of a (line, plane)-set which annihilates sets in $\phi^{-1}( \mathcal{C} \cup \mathcal{E} \cup \mathcal{F})$, a countable sum of quasi-invariant measures supported on (line, line)-sets or (plane, plane)-sets which annihilate sets from $\phi^{-1} (\mathcal{C})$ and a countable number of measures supported on (plane, line)-sets. The result then follows from Lemma \ref{ideal}, Proposition \ref{atom}, Lemma \ref{lem4.18} and Lemma \ref{lem4.19}. $\Box$

 \section{Statement of the result}
 We summarize the previous discussion in the next Theorem.
 Given two finite measures $\mu$, $\nu$ on a locally compact space $X$ we will write $\mu = \mu_{\ll \nu}  + \mu_{\perp \nu}$, where $\mu_{\ll \nu}$ is absolutely continuous with respect to $\nu$ and $\mu_{\perp \nu}$ is singular with respect to $\nu$.
 
 \begin{lem}
 \label{lem18}
 Let $X$ be a locally compact space endowed with an action of a discrete countable group $\Gamma$ and $\mu$ a finite Borel measure on $X$. There are a quasi-invariant finite measure $\tilde{\mu}$ on $X$ and a Borel set $A\subset X$ of full $\mu$-measure such that $\mu=\mu|_A = \tilde{\mu}|_A$. Moreover, $\tilde{\mu}$ is a Radon measure whenever $\mu$ is.
 \end{lem}
 \proof Fix an enumeration $\{\gamma_n\}_{n \in \mathbb{N}}$ of the group $\Gamma$. For every $n \in \mathbb{N}$ let $\mu_0 =\mu$, $\mu_n := (\gamma_n \mu)_{\perp \sum_{i=0}^{n-1} \mu_i}$ for $n>0$. Consider the measure $\tilde{\mu}=\sum_{n=0}^\infty 2^{-n}\mu_n$. We show that $\tilde{\mu}$ is quasi-invariant. Let $E$ be a Borel subset with $\tilde{\mu}(E)=0$. Then $\mu_n (E)=0$ for every $n$ and so $\sum_{i=0}^{n-1} \mu_i (E)=0$ for every $n$. Note now that if $\sum_{i=0}^{n-1} \mu_i (E)=0$ and $\mu_n (E)=0$ for some $n>0$ then $\gamma_n \mu (E)=0$, in fact
 \begin{equation*}
  (\gamma_n \mu) (E)=(\gamma_n \mu)_{\ll \sum_{i=0}^{n-1} \mu_i}(E) +( \gamma_n \mu)_{\perp \sum_{i=0}^{n-1} \mu_i}(E)=(\gamma_n \mu)_{\ll \sum_{i=0}^{n-1} \mu_i}(E) + \mu_n (E)=0.
 \end{equation*}
It follows that $(\gamma_n \mu) (E)=0$ for every $n$. Hence, given $\gamma \in \Gamma$, we have $\mu_n (\gamma E)\leq \mu (\gamma_n^{-1} \gamma E)=0$ for every $n \in \mathbb{N}$ and the result follows.\\
For every $0 < n \in \mathbb{N}$ there is a Borel set $A_n$ such that $\mu (X \backslash A_n)=0$, $\mu_n (A_n)=0$. Let $A = \bigcap_n A_n$ and let $E$ be a Borel set. Then $\tilde{\mu} (E \cap A) = \sum_{n=0}^\infty 2^{-n} \mu_n (E \cap A)=\mu (E \cap A)$. The regularity of the measure $\tilde{\mu}$ follows from the construction. $\Box$

 \begin{thm}
 \label{result}
 Every $\SL(3,\zz) \times \SL(3,\zz)$-quasi-invariant Radon measure $\mu$ on $\partial_\beta \SL(3,\zz)$ decomposes in the sum of two mutually singular quasi-invariant Radon measures $\mu_{\topp}$ and $\mu_{\lin}$ such that the associated Koopman representations satisfy:
 \begin{itemize}
 \item[-] $\pi_{\mu_{\topp}}$ is a $c'_0 (\SL(3,\zz) \times \SL(3,\zz))$-representation (see Definition \ref{def-irep}, Definition \ref{def-c0}),
 \item[-] $\pi_{\mu_{\lin}}$ is a $C^*_{red} (\SL(3,\zz) \times \SL(3,\zz))$-representation.
 \end{itemize}
In particular, the centralizer of every infinite subset of $\SL(3,\zz)$ is amenable and the same property does not hold for $\SL(n,\zz)$ if $n \geq 4$.
 \end{thm}
 \proof A quasi-invariant measure on $\partial_\beta \SL(3,\mathbb{Z})$ can be decomposed in the sum of a measure supported on $\phi^{-1} (M_{3,2})$, one supported on $\phi^{-1}(M'_{3,2})$ and one supported on $\phi^{-1} (Y_1)$. The part supported on $\phi^{-1} (Y_1)$ is tempered and so, in virtue of Proposition \ref{pinv-koop}, it is enough to study a measure supported on $\phi^{-1}(M_{3,2})$. Let $\mu$ be a quasi-invariant Radon measure on $\phi^{-1} (M_{3,2})$ and decompose it as a sum of a measure supported on preimages of (line-plane)-sets and a top-dimensional measure. The Koopman representations associated to the orbits of (line, plane)-sets are $C^*_{red}(\SL(3,\zz) \times \SL(3,\zz))$-representations in virtue of Proposition \ref{proposition4.20}, while the Koopman representation associated to the top-dimensional measure is a $c'_0 (\SL(3,\zz) \times \SL(3,\zz))$-representation in virtue of Proposition \ref{prop11}.\\
 Let now $\Gamma$ be any countable discrete group satisfying the property of the statement and let $\{x_n\}$ be an infinite subset of $\Gamma$; denote by $\Lambda$ the centralizer of $\{x_n\}$ in $\Gamma$. Let $\omega \in \partial_\beta \mathbb{N}$ and $\phi$ the associated state on $C(\partial_\beta \Gamma)$ given by $\phi (f) = \lim_{n \rightarrow \omega} \tilde{f}(x_n)$, where $\tilde{f}$ is any representative of $f$ in $C(\Delta_\beta \Gamma)$. Let $\mu_0$ be the associated Radon measure on $\partial_\beta \Gamma$ and $\mu$ its quasi-invariant extension as in Lemma \ref{lem18}. Suppose we have a decomposition $\mu=\mu_{\topp} + \mu_{\lin}$ as in the statement. If $\Lambda$ is finite, there is nothing to prove, hence we can suppose that $\Lambda$ is infinite. First observe that $\mu_0 \perp \mu_{\topp}$, indeed otherwise we could write $\mu_0 = \mu_1 + \mu_2$, where $\mu_1 \perp \mu_2$ and $ 0 \neq\mu_1 \leq \mu_{\topp}$. Since the action of $\Lambda$ on $L^2 (\mu_0)$ is trivial, the same is true for its restriction to $L^2 (\mu_1)$, but $L^2 (\mu_1)$ embeds $\Lambda$-equivariantly inside $L^2 (\mu_{\topp})$ and this is the desired contradiction. Hence we certainly have $\mu_0 \leq \mu_{\lin}$. Since the Koopman representation $\pi_{\mu_{\lin}}$ is a $C^*_{red} (\Lambda)$-representation, we have that the trivial representation of $\Lambda$ is a $C^*_{red}(\Lambda)$-representation, which entails amenability of $\Lambda$. As a consequence $\SL(n,\zz)$ does not satisfy this property for $n \geq 4$. $\Box$
 
As the anonymous referee pointed out, the fact that infinite subsets of $\SL(3,\mathbb{Z})$ have amenable centralizers can also be obtained by elementary algebraic methods (\cite{oza}, Example 7).
 
  \begin{cor}
 The Calkin representation $\pi_{\mathcal{C}}$ of $C^*(\SL(3,\zz) \times \SL(3,\zz))$ is a $C^*_{c'_0}(\SL(3,\zz) \times \SL(3,\zz))$-representation. Moreover, the same result does not apply to $\SL(n,\zz)$ for $n \geq 5$. 
 \end{cor}
 \proof The proof of Proposition \ref{proposition4.20} relies (see the proofs of Proposition \ref{s5}, Lemma \ref{lem4.18} and Lemma \ref{lem4.19}) on the fact that the global stabilizers of the preimages in $\partial_\beta G(\mathbb{Z})$ of sets from $\mathcal{D}$, $\mathcal{C}$, $\mathcal{E}$ and $\mathcal{F}$ act amenably on the respective preimages. Hence the same proofs apply for the respective preimages in $\psi_\omega^{-1}(\partial_\beta \SL(3,\mathbb{Z}))$ (see the discussion preceding Theorem \ref{thm2.4.3} for the definition of the map $\psi_\omega$). Therefore, every quasi-invariant finite Radon measure on $\psi_\omega^{-1}(\partial_\beta \SL(3,\mathbb{Z}))$ supported on the preimages of (line,plane)-sets is tempered. Also Proposition \ref{prop-inverse} admits a generalization to the case at hand, by considering the map ${}^* \inv : \sigma(l^\infty(\SL(3,\mathbb{Z})_\omega) \rightarrow \sigma(l^\infty(\SL(3,\mathbb{Z})_\omega)$ induced by $(f_n) \mapsto (f_n \circ \inv)$ and the first part of the proof of Theorem \ref{result} applies, showing that every $\SL(3,\mathbb{Z})\times \SL(3,\mathbb{Z})$-quasi-invariant finite Radon measure on $\psi_\omega^{-1}(\partial_\beta \SL(3,\mathbb{Z}))$ decomposes as the orthogonal sum of a tempered measure and a measure having  $c'_0$-matrix coefficients. The result follows from Lemma \ref{ideal} and Theorem \ref{thm2.4.3}.
 
 Let now $n \geq 5$ and suppose that the Calkin representation of $\SL(n,\zz) \times \SL(n,\zz)$ on $\mathbb{B}(l^2 (\SL(n,\zz))/\mathbb{K}(l^2 (\SL(n,\zz))$ is a $C^*_{c'_0} (\SL(n,\zz) \times \SL(n,\zz))$-representation. Embed $\SL(3,\zz) $ in the upper-left corner of $\SL(n,\zz) \times \SL(n,\zz)$ via the homomorphism 
 \begin{equation*}
 g \mapsto \left(\left( \begin{array}{cc}	g	&	0	\\
 				0	&	1	\end{array}\right), \left( \begin{array}{cc}	g	&	0	\\
 				0	&	1	\end{array}\right)\right).
				\end{equation*}
				In virtue of \cite{BrGu} Proposition 2.8 the restriction of the Calkin representation to $\SL(3,\zz)$ is a $C^*_{c_0} (\SL(3,\zz) ) $-representation. But this representation contains the trivial representation of $\SL(3,\zz)$ (c.f. \cite{Sk} Remark 4.6) and this would imply a-T-menability of $\SL(3,\zz)$. $\Box$

 \section{Acknowledgments}
The authors thank the anonymous referee for the comments on a previous version of this manuscript, in particular for suggesting a better proof of Lemma 2.5 and for bringing to their attention a simple, algebraic proof of the amenability of the centralizers of infinite sets in $\SL(3,\mathbb{Z})$. They also thank Prof. M. di Nasso and M. Pierobon for interesting discussions concerning Theorem \ref{thm2.4.3}. They acknowledge the support of INdAM-GNAMPA and the MIUR Excellence Department Project awarded to the Department of Mathematics, University of Rome Tor Vergata, CUP E83C180001000 and of the grant Beyond Borders: "A geometric approach to harmonic analysis and spectral theory on trees and graphs", CUP: E89C20000690005. The first named author was supported by the MIUR - Excellence Departments - grant: "$C^*$-algebras associated to $p$-adic groups, bi-exactness and topological dynamics", CUP: E83C18000100006 and by the grant Beyond Borders: "Interaction of Operator Algebras with Quantum Physics and Noncommutative Structure", CUP: E84I19002200005, during the period of this research. The second named author acknowledges the partial  support of the grant  "The convex space of sofic representations", CNCS Romania, PN-III-P1-1.1-TE-2019-0262. The present project is part of: - OAAMP - Algebre di operatori e applicazioni a strutture non commutative in matematica e fisica, CUP E81I18000070005. Florin R\u adulescu is a member of the Institute of Mathematics of the Romanian Academy.

\bibliographystyle{mscplain}
 \bibliography{biblio}

\baselineskip0pt

\end{document}